%%%% 
\input amstex

\def\b1{\text{\bf 1}}

\def\CA{{\Cal A}}

\def\CC{{\Cal C}}
\def\CD{{\Cal D}}
\def\CF{{\Cal F}}

\def\CH{{\Cal H}}

\def\CG{{\Cal G}}
\def\CL{{\Cal L}}

\def\CK{{\Cal K}}
\def\CO{{\Cal O}}

\def\CS{{\Cal S}}

\def\CV{{\Cal V}}

\def\CCm{{\Cal C^{\text{min}}}}

\def\#{\,\check{}}

\def\id{\text{id}}

\def\Ker{\text{Ker}}

\def\Spec{\text{Spec}}

%symbols

\def\hra{\hookrightarrow}
\def\iso{\buildrel\sim\over\rightarrow} 

\def\lra{\longrightarrow}

% space between paragraphs 
%\parskip=6pt

%\documentstyle{conm-p}
\documentstyle{amsppt}
\NoBlackBoxes

\topmatter\title    Constructible sheaves are holonomic       \endtitle \author A.~Beilinson \endauthor   

\dedicatory To  Joseph Bernstein
  \enddedicatory
 
\abstract We develop the theory of singular support for \'etale sheaves on algebraic varieties over an arbitrary base field.
\endabstract

\thanks The author was supported in part by NSF Grant DMS-1406734.\endthanks

 \subjclass Primary  14F05; Secondary 14C21 \endsubjclass


\address Department of Mathematics, University of Chicago, Chicago, IL 60637\endaddress \email sasha\@math.uchicago.edu\endemail

\endtopmatter

\document

In  \cite{KS} Kashiwara and Shapira, motivated by the theory of holonomic $\CD$-modules, have shown
that for any constructible sheaf $\CF$ on a    manifold $X$ 
each component of its singular support $SS(\CF )$ has dimension $\dim X$. Here
$SS(\CF )$ is
the smallest closed conical subset  of the cotangent bundle $T^* X$, such that  (locally on $X$) every  function $f$   with $df $ disjoint from $SS(\CF )$ is  locally acyclic relative to $\CF$. In this article we  establish a similar result for \'etale constructible sheaves on algebraic varieties over an arbitrary base field $k$. Our  tools are Brylinski's Radon transform \cite{B} (which is a global algebro-geometric version of the quantized canonical transformations from \cite{SKK}) and a version of the Lefschetz  pencils story \cite{Katz}.  

As was observed by Deligne \cite{D2}, if the characteristic of $k$ is finite then, in   contrast to the characteristic zero case  (either for constructible sheaves or   $\CD$-modules), the components of $SS(\CF )$ need not be Lagrangian  (see  examples in 1.5 below).

Recently Takeshi Saito \cite{S}    defined the characteristic cycle $CC (\CF )$, i.e.,  equipped the components of $SS(\CF )$ with multiplicities, and established its basic local and global properties   (the    Milnor type formula for the total dimension of   vanishing cycles  and the global Euler characteristic formula).

\enspace

 I am very grateful to Pierre Deligne, Volodya Drinfeld, Dennis Gaitsgory, Luc Illusie, Sasha Kuznetsov, and Takeshi Saito for valuable comments, suggestions, and  corrections.
 I learned about the singular support and the Radon transform from Joseph Bernstein at the end of the 1970s; this article is a small token of   gratitude.

\enspace

\head 1   Main results \endhead

\subhead{\rm 1.1}\endsubhead {\it Conventions and notation.}
 We fix a base field $k$; let $p$ be its characteristic. For us ``variety" means ``$k$-scheme  of finite type". A variety  $X$ is smooth if it is smooth relative to $k$; for a smooth $X$ we denote by $T^* X =T^* (X/k)$ the cotangent bundle relative to $k$. A map $f: X\to Y$ of smooth varieties yields a map of vector bundles $df: T^* Y_X := T^* Y\times_Y X \to T^* X$.
  For a vector bundle $\CV$ on $X$, we denote by $P(\CV )$ its projectivization (for $x\in X$ a point of the fiber $P(\CV )_x$ is a line in $\CV_x$). For a closed conical (i.e., $\Bbb G_m$-invariant) subset $C$ of $\CV$, we have its projectivization  $P(C)$ which is a closed subset of $P(\CV )$. The image of $C$ in $X$ is called the {\it base} of $C$; this is a closed subset of $X$.  A  {\it test pair}  on $X$  is a correspondence of type $X\buildrel{h}\over\longleftarrow U \buildrel{f}\over\lra Y$; we denote it by $(h,f) : X\leftarrow U \to Y$. 
 
 We fix   a  prime $\ell $ different from $p$;  ``sheaf" means ``bounded constructible complex of \'etale $\Bbb Z /\ell^n$-sheaves", $D(X)$ is the derived category of sheaves on a variety $X$.  A test pair $(h,f)$ as above is said to be {\it $\CF$-acyclic} for $\CF \in D(X)$ if $f$ is locally acyclic relative to $h^* \CF$ (see \cite{D1} 2.12). Thus $(\id_X ,\id_X )$ is $\CF$-acyclic if and only if $\CF$ is   locally constant, i.e.,  all the cohomology $\CH^i \CF$ are locally constant.
 
\subhead{\rm 1.2}\endsubhead 
 Let $X$ be a smooth variety, and let $C$ be a closed conical subset in $T^* X$.

  A map $h: U\to X$ with $U$ smooth is   {\it $C$-transversal} at a geometric point $u\in U$ if 
  $\Ker (dh_u )\cap C_{h(u)}\smallsetminus  \{0\}=\emptyset$.\footnote{The terminology of \cite{KS} 5.4.12 is  ``$h$ is non-characteristic for $C$."}
   A map   $f: X\to Y$ with $Y$ smooth is  {\it $C$-transversal } at a geometric point $x\in X$ if 
   $  (df_x)^{-1} (C_x )  \smallsetminus \{0\}  =\emptyset $.
We say that $f$ or $h$ is  $C$-transversal if it is $C$-transversal  at every geometric  point.

\proclaim{Lemma} (i) $C$-transversality is an open property:
for $f$  as above the set of  $x\in X$ such that $f$ is $C$-transversal at $x$ is open; ditto for $h$ and $u\in U$.
\newline
(ii) If $h$ as above is $C$-transversal, then the map 
$dh|_{C_U} : C_U := C\times_X U  \to T^* U$ is finite. Therefore its image $h^\circ C$ is a closed conical subset of $T^* U$.
\endproclaim

\demo{Proof} (i)  The preimage  $(df)^{-1}(C)$  is a closed conical subset of $T^* Y_X$. The image of   $P((df)^{-1}(C))$ is closed in $X$, and our set of $x$'s is its complement. In the case of  $h$, our set of $u$'s is the complement to the image of $P(C_U \cap \Ker(dh ))$ in $U$.  \newline
(ii) We can assume that $U$ is affine, $U=\Spec\, A$. Then $C_U =\Spec \, P$, $T^* U=\Spec\, R$  where
 $P$ and $R$ are $\Bbb Z$-graded $A$-algebras (due to the $\Bbb G_m$-actions), and $A=R/I$, $I:=R_{>0}$. The map    $dh|_{C_U}$ comes from a morphism of graded algebras $R\to P$. Since the gradings are nonnegative, a set of homogenous elements of $P$ generates $P$ as an $R$-module if (and only if) its image generates $P/IP$ as an $A$-module. Thus the map $dh|_{C_U}$  is finite if (and only if) it is finite over the zero section of $T^* U$, which is the $C$-transversality condition. \hfill$\square$ \enddemo

 Notice that if $h$ is smooth, then it is automatically $C$-transversal   and $h^\circ C =C_U$. If $f$ is $C$-transversal, then it is smooth on a neighborhood of the base of $C$.

\enspace

A  test pair $( h,f)$ as in 1.1  is said to be {\it $C$-transversal} at a geometric point $u\in U$  if $U$ and $Y$ are smooth, $h$ is $C$-transversal at $u$ and $f$ is $h^\circ C_u$-transversal at $u$. By the lemma, $C$-transversality is an open property with respect to $u\in U$. We say that
 $(h,f)$ is $C$-transversal at a subset $S$ of $U$ if it is $C$-transversal  at every $u\in S$; if $S=U$, then we simply say that  $(h,f)$ is $C$-transversal (which means that $h$ is $C$-transversal and $f$ is $h^\circ C$-transversal). Notice that a $C$-transversal test pair  is $C'$-transversal for every $C'\subset C$.
 
\remark{Examples} (i) If $C$ is the zero section of the cotangent bundle of $X$, then  $(h,f)$ is $C$-transversal if and only if $f$ is smooth. \newline (ii)   $(h,f)$ is $T^* X$-transversal if and only if the map $h\times f: U\to X\times Y$ is smooth. 
 \endremark

\enspace

For a  map of smooth varieties
$r: X\to Z$ and a closed conical subset $C$ of $T^* X$ whose  base   is proper over $Z$,  we denote by $r_\circ C$ the closed conical subset of $T^* Z$ defined as the image of $(dr)^{-1}(C)\subset T^* Z\times_Z X$ by the projection $T^* Z \times_Z X \to T^* Z$.

\subhead{\rm 1.3}\endsubhead  Our aim is to assign to each sheaf $\CF$ on a smooth variety $X$ a closed conical subset $SS(\CF )$ of $ T^* X$, so that for every $SS(\CF )$-transversal $h:U\to X$ one has $SS(h^* \CF )\subset h^\circ SS(\CF )$,  every $SS(\CF )$-transversal $f:X\to Y$  is locally acyclic relative to $\CF$, and to do this in the most economic way possible. We proceed as follows:

Let $X$ be a smooth variety, $\CF\in D(X)$.  We say that  $\CF  $ is {\it   micro-supported} on a closed conical subset $C$ of $ T^* X$   if  every  $C$-transversal  test pair     is $\CF$-acyclic. 
  
Consider the set $ \CC (\CF )$ of all $C$'s such that $\CF$ is micro-supported on $C$. If $C\in  \CC (\CF )$ and $C'\supset C$, then   $C'\in \CC (\CF)$. Our $ \CC (\CF )$ is not empty:

\proclaim{Lemma}
One has $T^* X\in \CC (\CF )$.
\endproclaim

\demo{Proof} Let $p_X, p_Y : X\times Y \to X, Y$ be the projections. Then
$p_Y$ is locally acyclic relative to   $p_X^* \CF$ by \cite{D1} 2.16. We are done, since for any test pair $(h,f)$ as in 1.1 
 one has $h^* \CF = (h\times f)^* (p_X^* \CF )$, and  for a $T^* X$-transversal  $(h,f)$,  the map $h\times f$ is smooth (see Example (ii) in 1.2). \hfill$\square$
\enddemo

 We say that $\CF$ {\it has singular support} if  $ \CC (\CF )$  has the smallest element. The latter is 
 denoted then
 by $SS(\CF )=SS(\CF,X/k)$ and called the {\it singular support}, or {\it micro-support} of $\CF$.

\proclaim{ Theorem} (i) Every $\CF$ has singular support. 
\newline  
(ii) For a connected $X$ each irreducible component of $SS(\CF )$ has dimension $\dim X$. \endproclaim

For the proof of (ii) see 4.10; (i) and the
 upper bound for $\dim SS(\CF )$ are  in \S 3.

\remark{Remark} If the characteristic $p$ of $k$ is zero then, by Kashiwara-Schapira \cite{KS}, every component of $SS(\Cal F)$ is a Lagrangian cone.\footnote{Proof (Kuznetsov). (a) By Theorem 1.3(ii) it is enough to  show that
$SS(\Cal F)$ is isotropic. 
(b) By   1.4(ii) we can assume that $\CF$ is an irreducible perverse sheaf, hence it is the middle extension of a local system $\CF_\eta$ at the generic point of a closed irreducible subscheme $Y$ of $X$. (c) By   1.4(ii) it is enough to find some $\CF' \in D(X)$ supported on $Y$ such that $\CF'_\eta$ is a local system that contains $\CF_\eta$ and $SS(\CF')$ is isotropic. (d) To that end we find, by de Jong,   a proper map $r: Z\to X$ with $Z$   smooth and irreducible of dimension $\dim Y$ such that  $r(Z)=Y$ and  the local system $r^* \CF_\eta$ on the generic point of $Z$ is trivial.  Let $\CG$ be the constant sheaf on $Z$ with fiber $r^* \CF_\eta$. Then
$\CF':=r_* \CG$ is the sheaf promised in (c): indeed, 
$SS(\CG )$ is the zero section of $T^* Z$ by 2.1(iii) hence is isotropic,
$SS(\CG' )\subset r_\circ SS(\CG )$ by 2.2(ii), and $r_\circ$
 sends isotropic cones to isotropic cones.}
\endremark

\remark{Question {\rm (Drinfeld)}}        Which conical irreducible subsets 
$C\subset T^* X$ of dimension $\dim X$ can be realized as   components
of $SS(\CF )$ for some $\CF \in D(X)$ in case $p>0$?
\endremark

\enspace

By  Deligne \cite{D3}, if $X$ is a surface, then {\it every}  $C$ can be realized in this manner.\footnote{Deligne shows that every $C$ that is not a conormal can be identified \'etale locally at the generic point with one of $C_n$'s from the next example.}

\remark{Example} Set $V:= \Spec\, k [ x,t]$; for $g\in k[x,t]$ let $r_g : V\to V$ be
 the map $(x,t)\mapsto (g(x,t),t)$. Consider the polynomials $g_0 (x,t):= x^p +x t^2$ and $g_n (x,t) := x^{p^n} +xt$, $n\ge 1$; set $\CF_n := r_{n *} \Bbb Z/\ell _V \in D(V)$ where $r_n := r_{g_n}$. Then $SS(\CF_n )$ is the union of
the zero section of $T^* V$ and the cone $C_n$ over the $x$-axis $(t=0)$ generated by the section $dx$ if $n=0$ and $dx- x^{1/p^n}dt$ if $n\ge 1$.\footnote{Proof.  The map $r_n$ is finite,  so $SS(\CF_n )\subset r_{n\circ} SS (\Bbb Z/\ell _V)$  by 2.2(ii). Since $SS(\Bbb Z/\ell _V)$ is the zero section of $T^* V$ (see 2.1(iii)) and $r_n$ is \'etale over the complement to the $x$-axis, an immediate computation shows that $r_{n\circ} SS (\Bbb Z/\ell _V )$ is the union of $C_n$ and the zero section of $T^* V$.  Now $SS(\CF_n )$ contains the zero section of $T^* V$ since $\CF_n$ is nonzero at the generic point of $V$, and $SS(\CF_n )$ is not equal to the zero section at the generic point $\eta_x$ of the $x$-axis since $\CF_n$ is not locally constant there (use 2.1(iii)). We are done since every closed subcone of $C_n$ other than $C_n$ is contained in the zero section at $\eta_x$.
} Notice that $C_0$ is not Lagrangian and $C_n$, $n\ge 1$, is a Lagrangian cone that  is not   conormal to its base. \endremark

\remark{Exercise} An irreducible Lagrangian cone $C\subset T^* X$ with smooth base $Y$ coincides with the conormal bundle to $Y$ if (and only if) the projection $C\to Y$ is smooth at the generic point of $C$.
\endremark

\subhead{\rm 1.4}\endsubhead  Here are some other useful properties of $SS$ (to be proved in 4.10):

\proclaim{ Theorem} (i) For a smooth map $\pi : Z\to X$ one has $SS (\pi^* \CF  )=\pi^\circ SS(\CF )$.  
\newline (ii) If  $\{ \CF_\alpha \}$ are perverse Jordan--H\"older constituents\footnote{I.e., $\{ \CF_\alpha \}$ are irreducible perverse sheaves that can be realized as subquotients of some perverse sheaf cohomology ${}^p \CH^a \CF$.}  of $\CF$, then $SS(\CF )\! =\! \cup_\alpha SS (\CF_\alpha )$.  
(iii) Compatibility with the base field change: Let $k'/k$ be any extension of the base field,  and let $\CF_{k'}$ be the pullback of $\CF$ to $X_{k'}:=X\otimes_k k'$. Then
 the closed subsets $SS(\CF, X/k )_{k'}$ and $SS(\CF_{k'}, X_{k'}/k' )$ of  $(T^* X)_{k'}=T^* ( X_{k'}/k')$ coincide.
\endproclaim

\subhead{\rm 1.5}\endsubhead  The next apparent weakening of the above notion is useful. We say that $\CF$ is {\it weakly micro-supported} on a closed conical subset  $C$ of $T^* X$ if every $C$-transversal test pair $(h,f)$  {\it that satisfies the next two extra conditions} is $\CF$-acyclic:
\newline
- $f$ is a function, i.e., $f:U\to Y=\Bbb A^1$; \newline
- if $k$ is infinite, then $h:U\to  X$ is an open embedding; if $k$ is finite, then $h$ is the composition $U=V_{k'}:=V\otimes_k k' \to V \hra X$ where $V$ is an open subset of $X$ and $k'$ is a finite extension of $k$.

The set $\CC' (\CF )$ of all closed conical subsets of $T^* X$  on which $\CF$ is weakly micro-supported is   a filter\footnote{To check that  $C_1 ,C_2 \in \CC' (\CF)$ implies $C_1 \cap C_2  \in \CC' (\CF)$ notice that if a test pair $(h,f)$ as in 1.1 with $\dim Y=1$ is $C_1 \cap C_2$-transversal, then locally on $U$ it is either $C_1$- or $C_2$-transversal. The latter assertion need not be true if $\dim Y>1$ (consider cones
$C_1$, $C_2$ that are nonzero with zero intersection at the generic point of $X$ and
 $(h,f)=(\id_X ,\id_X )$).} (the   assertion with ``weakly micro-supported" replaced by ``micro-supported" is 1.3(i), and it is not evident). Denote by $SS^w(\CF)$ its minimal element.  Explicitly, $SS^w(\CF )$ is the closure in $T^* X$ of the set of all points $(x, df(x))$ where $x$ is a closed point of $X$ and $f$ is a function on a Zariski neighborhood\footnote{With modification as above in case $k$ is finite. Indeed, the modification is needed to ensure that $SS^w (\CF )$ for $\CF$ a skyscraper sheaf at $x\in X$ equals $T^*_x X$.} of $x$ which is not locally acyclic relative to $\CF$ at $x$. 

Clearly $SS^w (\CF )\subset SS(\CF )$ so  1.3(ii) implies that  $\dim SS^w (\CF )\le \dim X$; this assertion was conjectured by Deligne in \cite{D2}. The next result is proved in 4.9:

\proclaim{ Theorem}  One has $SS^w(\CF )=SS(\CF )$. \endproclaim

\subhead {\rm 1.6}\endsubhead  {\it The Radon and Legendre transforms: a reminder.}  The main tool used in the proofs is Brylinski's geometric Radon transform. We recall the definition and list the key properties that we will need.   See \cite{B} for details. 

Let $V$ be a vector space of dimension $n+1$, $V^\vee$ its dual. Let $\Bbb P$, $\Bbb P^\vee$ be the corresponding projective spaces, let $Q\hra \Bbb P \times \Bbb P^\vee $ be the incidence correspondence, and let $p, p^\vee : Q \rightrightarrows \Bbb P,\,\Bbb P^\vee \!$ be the projections.  For $x\in \Bbb P$, $x^\vee \in \Bbb P^\vee$ we denote by $Q_x$, $Q_{x^\vee}$ the corresponding fibers (so $x$ is a line in $V$ and $Q_{x}$ is the hyperplane in $\Bbb P^\vee$ formed by lines in $V^\vee$ orthogonal to $x$). For $(x,x^\vee )\in Q$ the tangent spaces to the fibers $Q_x$, $Q_{x^\vee}$
at $(x,x^\vee )$  intersect trivially, so the orthogonal complement to $T_{(x,x^\vee )}Q_x \oplus T_{(x,x^\vee )}Q_{x^\vee} \subset T_{(x,x^\vee )}Q$ is a line $\lambda_{(x,x^\vee )}\subset T^*_{(x,x^\vee )}Q$.

\subhead{\rm 1.6.1}\endsubhead We have the
{\it  Radon transform} functors $R:= p_{ *}^\vee  p^* [n-1] : D(\Bbb P )\to D(\Bbb P^\vee )$ and $R^\vee := p_* p^{\vee *}[n-1] : D(\Bbb P^\vee )\to D(\Bbb P )$. For  $\CK \in D(\Bbb P^\vee )$ one has: \newline (i) $R$ and $R^\vee (n-1)$ (the Tate twist) are (both left and right) adjoint; the cones of the adjunction  maps $ RR^\vee (\CK )(n-1) \to \CK$, $\CK \to R^\vee R  (\CK )(n-1) $ are locally constant.\footnote{Since we live on a projective space this means that the cohomology sheaves come from Spec$\, k$.} \newline (ii) If $\CK$ is a perverse sheaf,   then the perverse cohomology sheaves ${}^p \CH^i R^\vee (\CK )$ are locally constant for $i\neq 0$. If $\CK$ is locally constant, then so is $R^\vee (\CK )$.
\newline (iii) If $\CK$ is an irreducible perverse sheaf which is not locally constant, then
${}^p \CH^0 R^\vee (\CK )$ has single not locally constant irreducible perverse sheaf in its Jordan--H\"older series.

\subhead{\rm 1.6.2}\endsubhead 
 One has the   {\it Legendre transform}  identifications  
$ P(T^* \Bbb P ) \buildrel{\sim}\over\leftarrow Q\iso P(T^* \Bbb P^\vee )   $
compatible with the projections to $\Bbb P$ and $\Bbb P^\vee$ respectively. 
The left arrow assigns to $(x,x^\vee )\in Q$ the point $(\lambda_{{x^\vee} x})\in P(T^*  \Bbb P )_x$ that corresponds to the 
conormal line
$\lambda_{{x^\vee} x}\subset T^*_x  \Bbb P$   to the hyperplane $Q_{x^\vee}$ at $x$.  The right arrow is  defined dually.

  For $(x,x^\vee )\in Q$ consider the embeddings $ T^*_x \Bbb P \buildrel{dp}\over\lra   T^*_{(x,x^\vee )} Q \buildrel{dp^\vee}\over\longleftarrow T^*_{x^\vee} \Bbb P^\vee .  $  The two vector subspaces intersect transversally by the line $\lambda_{(x,x^\vee )}$, and  one has $dp (\lambda_{{x^\vee} x})=\lambda_{(x,x^\vee )}=dp^\vee (\lambda_{{xx^\vee} })$. Thus $T^*_{(x,x^\vee )} Q$ is the direct sum of $
T^*_x \Bbb P$ and $ T^*_{x^\vee} \Bbb P^\vee$ with the lines $ \lambda_{{x^\vee} x}$ and $\lambda_{{xx^\vee} }$ identified.

\remark{Remark} The line subbundle $\lambda $ of $ T^* Q$ is a contact structure on $Q$, and the Legendre transforms identify it with the canonical contact structures on the projectivizations of the cotangent bundles to  $\Bbb P$ and $\Bbb P^\vee$.
\endremark

\subhead{\rm 1.7}\endsubhead Let us explain   how to recover, after a  Veronese embedding,  the singular support of any sheaf on $\Bbb P$  from  the ramification divisor on the Radon transform side. The theorems in 1.3, 1.4, and 1.5 are deduced from this.

We return to the setting of 1.6. 
 Let $i=i_d : \Bbb P\hra \tilde{\Bbb P} $ be the Veronese embedding of degree $d\ge 2$. Thus $\tilde{\Bbb P}$ is a projective space of
 dimension $N_d =  {n+d \choose d} -1$ and points of $\tilde{\Bbb P}^{\vee}$ are degree $d$ hypersurfaces in $\Bbb P $. Let $\tilde{p}, \tilde{p}^\vee : \tilde{Q} \rightrightarrows \tilde{\Bbb P},\tilde{\Bbb P}^{\vee}$ be the incidence correspondence and let $\tilde{R}$, $\tilde{R}^\vee$ be the Radon transforms. 
  
A closed conical subset $C$ of $T^* \Bbb P$ yields the closed conical subset $i_\circ C$ of $T^* \tilde{\Bbb P}$; let $P(i_\circ C)\subset P( T^* \tilde{\Bbb P})=\tilde{Q}$ be its projectivization, and let $D_{C}$ be the image of  $P(i_\circ C)$ in $\tilde{\Bbb P}^\vee$. 

Let $\CF\in D(\Bbb P )$ be any sheaf. Consider the Radon transform   $\tilde{R} (i_* \CF )\in D(\tilde{\Bbb P}^\vee )$. Let $D_\CF$ be the smallest closed subset of $\tilde{\Bbb P}^\vee $ such that 
 $\tilde{R} (i_* \CF )$ is locally constant on the complement $\tilde{\Bbb P}^\vee \smallsetminus D_\CF .$  
  
 \proclaim{Theorem}  $D_\CF$ is a divisor. For every its irreducible component $D_\gamma$ there is a unique irreducible closed conical subset $C_\gamma $ in $ T^* \Bbb P$ of dimension $n$ with $D_\gamma = D_{C_\gamma }$. 
  One has $SS(\CF )=\cup_\gamma \, C_\gamma$.  The maps $\tilde{p}^\vee_\gamma :  P(i_\circ C_\gamma )\to D_\gamma$ are generically radicial.  For $k$ perfect
 $\tilde{p}^\vee_\gamma$ are birational unless $p=2$ when the generic degree can also be 2.
  \endproclaim 

The proof of the theorem minus the last assertion is in 4.6 (case $d\ge 3$) and 4.7 (case $d=2$);  the last assertion (in a more precise form) is checked in 4.13.

\remark{Questions} (i) (Illusie) What would be the Picard--Lefschetz formula --- how to describe  $\tilde{R} (i_* \CF )$ 
on an \'etale neighborhood of  the generic point of $D_\alpha$?  \newline (ii) More generally, if $U$ is the strict Henselization of a variety at the generic point of a divisor $D$, then 
can one describe explicitly the category  of those sheaves on $U$ whose singular support is the union of the zero section and a given  $C\subset T^* U|_{D_U}$ with $P(C)$   radicial over $D_U$? 
\newline (iii) What would be the theory of  microlocal perverse sheaves in the present setting? Does the codimension-three conjecture hold? \endremark

\head 2      Some  elementary  lemmas \endhead  

We collect several  simple facts to be used in the proofs of the theorems from \S 1. For this section   $X$ is a smooth variety, $C $ is a closed conical subset of $T^* X$, and $\CF$ is a sheaf on $X$.
 
\proclaim{{\rm 2.1.} Lemma} (i) The base of $SS^w (\CF )$ (see 1.5)   equals the support of $\CF$.
\newline
(ii) If $C\in\CC (\CF)$,  then 
for every $C$-transversal test pair $(h,f): X\leftarrow U\to Y$  the map $f$ is   universally locally acyclic relative to $h^* \CF$. \newline (iii) 
$\CF$ is  micro-supported on the zero section of $T^* X$ if and only if $\CF$ is  locally constant.\newline (iv)  All sheaves that are micro-supported on $C$ form a thick subcategory of $D(X)$.
\endproclaim

\demo{Proof} (i) The support $Z$ of $\CF$ evidently contains the base $B$ of
$SS^w(\CF )$. Since $B$ is a closed subset (see 1.1)  it is enough, replacing $X$ by $X\smallsetminus B$, to show that for $\CF \neq 0$, one has $SS^w(\CF )\neq \emptyset$, which is clear since the test pair $(h,f)=(\id_X , 0 )$ as in 1.5  is not $\CF$-acyclic.

(ii)  We need to check that  local acyclicity of $f$ relative to $h^* \CF$ remains true after any base change $g: Z\to Y$. First notice that it is enough to consider the case of smooth $g$: Indeed, we can assume that $Y$ and $Z$ are affine and, 
  by a standard argument, that the scheme $Z$   has finite type over $k$. Then $g$ 
 can be presented  as a composition $Z\hra Z':= \Bbb A^N_Y \to Y$ with the first map being a closed embedding and the second one  the projection. 
Then $U_Z :=U\times_Y Z$ is the preimage of $Z$ in $U_{Z'}:=U\times_{Y}Z'$, so the local acyclicity of $U_Z \to Z$ relative to the pullback of $\CF$
follows from that  of $U_{Z'}\to Z'$, and we are done.

 Since $C\supset SS^w(\CF )$ we can, by (i), 
  replace $U$ by a neighborhood of the support of $h^* \CF$ where $f$ is smooth. Now if $Z/Y$ is smooth,   then  $X \leftarrow U_Z \to Z$ is a $C$-transversal pair, so $U_Z \to Z$ is locally acyclic relative to the pullback of $\CF$. 
  
 (iii) If $\CF$ is  micro-supported on the zero section of $T^* X$,  then the test pair $(\id_X ,\id_X )$ is $\CF$-acyclic so $\CF$ is locally constant (see 1.1). If $\CF$ is locally constant,  then it is micro-supported on the zero section due to local acyclicity property of smooth maps (see Example (i) in 1.2).
 
 (iv) We need to check that our full subcategory is triangulated and idempotently closed. Both assertions follow directly from the definition.
 \hfill$\square$
\enddemo

 \subhead{\rm 2.2}\endsubhead  
For the definitions of $\pi^\circ C$ and $r_\circ C$ see 1.2.
 
 \proclaim{  Lemma}  
 (i)  If $\phi: Q\to X$ is $C$-transversal  and   $C\in\CC (\CF )$, then $ \phi^\circ C \in \CC(  \phi^* \CF)$. \newline
(ii)  Suppose we have $r: X\to Z$ with $Z$ smooth and  $C\in\CC (\CF )$ with base  proper over $Z$. Then $r_\circ C \in\CC (r_* \CF)$.   \endproclaim
 
\demo{ Proof} (i) is evident.  (ii) We need to check that any $r_\circ C$-transversal test pair  $(h', f'):  Z \leftarrow V \to Y$  is $r_* \CF$-acyclic. Let $r' : X_V :=X\times_Z V \to V$ and $h: X_V \to X$ be the projections.    Our $C$ contains by 2.1(i)  (since $C\supset SS^w(\CF)$) the image of the support of $\CF$ by the zero section of $T^* X$. Therefore the $r_\circ C$-transversality of $h'$ implies that the variety $X_V$ is smooth on some neighborhood $U$ of the support of $h^* \CF$. The test pair $(h|_U , f'r'|_U ): X\leftarrow U \to Y$ is $C$-transversal so $f'r'|_U$, hence $f'r'$, is locally acyclic relative to  $h^* \CF$. So, since $r'$ is proper on the support of $h^* \CF$, $f'$ is locally acyclic relative to $r'_* h^* \CF$ (see   3.13(i) in \S 3 below for the details). We are done since $h^{\prime *}r_* \CF =r'_* h^* \CF$ by proper base change.  \hfill$\square$
\enddemo

 \subhead{\rm 2.3}\endsubhead  Let $\CCm (\CF )$ be the subset of minimal elements of $\CC (\CF )$ (see 1.3).\footnote{Theorem 1.3(i) asserts that   $\CCm (\CF )$ has a single element. But we did not prove it yet.} 
 Since $T^* X$ is Noetherian, every $C\in \CC (\CF )$ contains some $C' \in \CCm (\CF )$. 
  
 \proclaim{Lemma} (i) Let $V$ be an open subset of $X$. If 
 $C\in \CCm (\CF)$, then $C_V\in \CCm (\CF|_V)$. The map $\CCm (\CF)\to 
\CCm (\CF|_V)$, $C\mapsto C_V$, is surjective.
\newline (ii) If  $\pi : Q \to X$  is   smooth surjective and $\pi^\circ C \in \CCm (\pi^* \CF)$,   then $C\in\CCm (\CF)$.
\newline (iii) If $C\in\CCm (\CF)$,   then the base of $C$ equals the support of $\CF$.
\endproclaim
  
\demo{Proof}  (i) Use  the next fact which follows from 1.2(i): If $\CF$ is micro-supported on $C$ and $\CF|_V$ is micro-supported on $C'_V$,  then $\CF$ is micro-supported on the union of $C_{X\smallsetminus V}$ and the closure of $C'_V$ in $T^* X$. (ii) follows from 2.2(i).  To prove (iii) notice that  the base of $C$  lies in the support of $\CF$ by (i), and we are done by 2.1(i) since $C\supset SS^w (\CF)$.
\hfill$\square$
\enddemo

 \subhead{\rm 2.4}\endsubhead A closed conical subset $C$ of a vector bundle $\CV $ over a variety $S$ is  {\it strict} if none of its irreducible components   lies in the zero section $S\hra \CV$. The map $C\mapsto P(C)$   is a bijective correspondence between strict closed conical subsets of $\CV$ and closed subsets of $P(\CV )$; the inverse map sends a closed subset of $P(\CV )$ to the cone over it.  
 
\proclaim{Lemma} Suppose  $X$ is connected and $C\in\CCm (\CF)$. If $\CF$ vanishes at the generic point of $X$ then $C$ is strict; otherwise $C$ is the union of a strict subset (which is the cone over $P(C)$) and the zero section of $T^* X$. 
\endproclaim

\demo{Proof} By 2.3(i) we can replace $X$ by the complement to the image of $P(C)$. Then  $\CF$ is locally constant by 2.1(iii) and we are done by 2.3(iii). \hfill$\square$
\enddemo

\subhead{\rm 2.5}\endsubhead Let $i: X\hra P$ be a closed embedding of smooth varieties.

\proclaim{Lemma} (i) One has $SS^w(i_* \CF )= i_\circ SS^w(\CF )$. \newline
(ii) If the sheaf $i_* \CF$ on $P$ has singular support  (see 1.3), then $\CF$ has singular support and,
if $k$ is infinite, one has   $SS(i_* \CF )= i_\circ SS(\CF )$. If $k$ is finite, then
the latter assertion is true if we know that for every finite extension $k'/k$ the sheaf $i_{k'*} \CF_{k'}$ on $P_{k'}$ has singular support. 
\endproclaim

\demo{Proof}  (i) follows since (a) a function $g$ on $P$ is locally acyclic relative to $i_* \CF$ if and only if $g|_X$  is locally acyclic relative to $ \CF$, and (b) for every datum $(f,x,\nu)$, where $f\in \CO (X)$, $x\in X$ is a closed point, and $\nu\in T^*_{x}P$ is such that $di (\nu )\in T^*_x X$ equals $ df_x $, one can find Zariski locally  a function $g$ on $P$ such that $g|_X =f$ and $dg_{x}=\nu$.  

(ii) By 2.3(i) the claim is Zariski local. So we can assume that $P$ is affine and  there is an \'etale map $\chi =(\chi_n ,\chi_{m-n}): P\to\Bbb A^m =\Bbb A^n \times \Bbb A^{m-n}$ with $X=\chi_n^{-1}(\Bbb A^n)$. Our $\chi$  yields a splitting $s=s_\chi : T^* X \hra T^* P|_X$ of  $di :T^* P|_X \twoheadrightarrow T^* X$. 

It also provides a datum $(\tilde{P},\phi,\rho )= (\tilde{P}_\chi,\phi_\chi ,\rho_\chi )  $ where $\phi : \tilde{P}\to P$
is an \'etale map  with $\phi^{-1}(X)\iso X$ (so we have $X\hra \tilde{P}$ that lifts $ i$) 
and $\rho : \tilde{P} \to X$ is  
 a retraction  such that $d\rho|_X : T^* X \to T^* P|_X $ equals $s$. 
 Namely,  consider $\tilde{P}' :=P\times_{\Bbb A^n}X$ where the maps to $\Bbb A^n$ are $\chi_n$ and $\chi_n |_X$; let $\phi' : \tilde{P}' \to P$ and $\rho' : \tilde{P}' \to X$ be the projections. Then $\phi^{\prime -1}(X)=X\times_{\Bbb A^n}X$ is the disjoint union of the diagonal $X$ and its complement $K$. Our $\tilde{P}$ is $\tilde{P}'\smallsetminus K$, $\phi$ and $\rho$ are the restrictions of $\phi'$ and $\rho'$ to $\tilde{P}$.

(a)  Let us show that  $\CF $ has singular support.  By 2.3(iii) the base of $SS(i_* \CF )$ lies in $X$ so we have a closed conical subset $C:=\rho_\circ \phi^\circ SS(i_* \CF )= s^{-1} (SS(i_* \CF ))$ of $T^* X$. Then $SS (\CF )$ equals $C$: Indeed, $C\in \CC (\CF )$ by 2.2 since $\CF = \rho_* \phi^* (i_* \CF)$, and it is the smallest element of $\CC (\CF )$ since for any $C'\in \CC (\CF )$ one has $i_\circ C' \in \CC (i_* \CF )$ by 2.2(ii), so $i_\circ C' \supset SS(i_* \CF )$ and  $C' =\rho_\circ \phi^\circ (i_\circ C')\supset C$.

(b) Let us prove that  $SS(i_* \CF )= i_\circ SS(\CF )$. By 2.2(ii) one has $SS(i_* \CF )\subset  i_\circ SS(\CF )$. By (i) and 2.1(i) our $SS(i_* \CF )$ contains  the conormal bundle $T^*_{X}P$ restricted to the support of $\CF$ (recall that $SS(i_* \CF )$ contains $SS^w (i_* \CF )$). Suppose $V:= i_\circ SS(\CF )\smallsetminus SS(i_* \CF )\neq \emptyset$. Then we can find $x\in X$ and $\mu \in SS(\CF )_x \smallsetminus \{0\}$  such that the preimage in $V_x$ of $\mu$ by $  T^*_x P  \buildrel{di}\over\lra T^*_x X$ is not empty.

Suppose   $k$ is infinite. 
Then, replacing $\chi$ by $g\chi$ for a sufficiently general matrix $g$ of type $g=\id_{k^m} +A \in GL_m (k)\subset \text{Aut}(\Bbb A^m )$, $A\in \text{Mat}_{n,m-n}(k) \subset \text{Mat}_{m,m}(k)$, we can ensure\footnote{Use the fact that for any   $\nu\neq 0$ in $\Bbb A^n$, the set $\{ A(\nu), A\in \text{Mat}_{n,m-n}(k)\}$ is dense in $\Bbb A^{m-n}$.} that
 $\mu \notin s^{-1}( SS (i_* \CF )_x )$ which is  $SS (\CF )_x$ by (a); contradiction.

Suppose $k$ is finite. Consider the base change of our picture to a finite extension  $k'/k$ of degree
prime to $\ell$, so we have the projection $\pi_X : X_{k'}\to X$, etc. 
Then $\CF$ is naturally a direct summand of $\pi_{X*} \CF_{k'}=\pi_{X*} \pi_X^* \CF$, so $SS(\CF )\subset SS(\pi_{X*} \CF_{k'}) \subset \pi_{X\circ} SS(\CF_{k'})\subset \pi_{X\circ} \pi_X^\circ SS(\CF)=SS(\CF)$ (see 2.2).
Thus $SS(\CF )=\pi_{X\circ} SS(\CF_{k'})$; since $SS(\CF_{k'})$ is $\text{Gal} (k'/k)$-invariant, one has $SS(\CF_{k'})=SS(\CF )_{k'}$. Ditto for sheaves on $P$, so $SS(i_{k' *}\CF_{k'}) =SS(i_* \CF)_{k'} $.
An immediate version of the above argument for infinite $k$ shows then that for large enough $k'$ one has $SS(i_* \CF )_{k'}=i_\circ SS(\CF )_{k'}$ hence
$SS(i_* \CF )=i_\circ SS(\CF )$. \hfill$\square$\enddemo

\head 3  The upper estimate for $\dim SS(\CF )$ \endhead

\subhead{\rm 3.1}\endsubhead In this section we  prove the next weak version of the theorem in 1.3:

\proclaim{Theorem} 
Every sheaf  $\CF$ on a smooth variety $X$ has singular support. One has $\dim SS(\CF ) \le\dim X$. 
\endproclaim

\demo{Proof}  First notice that it is enough to consider the case when  $X$ is a projective space $\Bbb P=\Bbb P^n$. Indeed, the theorem for $\Bbb P^n$ implies that for $\Bbb A^n $ by 2.3(i), then for a smooth closed subvariety of  $\Bbb A^n $ by 2.5(ii), then for an arbitrary $X$ by 2.3(i).

Now we live on a projective space $\Bbb P$. The remaining proof has two steps:  first  we give a concrete description of $SS(\CF)$, then, using it, prove the estimate.
\enddemo

\subhead{\rm 3.2}\endsubhead    For a map $g: Y\to Z$ and a sheaf $\CG$ on $Y$ we denote by $E_g (\CG )$ the smallest closed subset of $Y$ such that $g$ is	 universally locally acyclic relative to $\CG$ on $Y\smallsetminus E_g (\CG )$. It has \'etale local nature with respect to $Y$ and $Z$.

 \proclaim{Theorem} Every  $\CF\in D(\Bbb P  )$ has singular support. Its projectivization $P(SS (\CF )) \subset  P(T^* \Bbb P )$ equals the Legendre transform of   $E_{p } (p^{\vee *}R  (\CF )  )\subset Q$. If $\CF$ vanishes at the generic point of $\Bbb P $, then $SS(\CF )$ is the cone over $P(SS (\CF ))$; otherwise $SS(\CF )$ is the union of this cone and  the zero section of $T^* \Bbb P $. \endproclaim

The proof is in 3.5;  it is based on general lemmas 3.3 and 3.4.

 \subhead{\rm 3.3}\endsubhead   The Legendre transform 1.6.2 yields a bijective correspondence between closed subsets of $P(T^* \Bbb P )$ and $P(T^* \Bbb P^\vee )$. For a closed conical  $C\subset T^* \Bbb P$ we denote by $C^\vee \subset T^* \Bbb P^\vee$  the cone over the Legendre transform of $P(C)\subset P(T^* \Bbb P )$, and by $C^+ \subset T^* \Bbb P$   the union of $C$ and the zero section of $T^* \Bbb P$. We use the same notation for $\Bbb P$ interchanged with $\Bbb P^\vee$; thus $C = C^{\vee\vee}$ if $C$ is strict (see 2.4) and $C^{\vee\vee +}=C^+$.

\proclaim{Lemma}  A sheaf $\CF \in D(\Bbb P  )$ is micro-supported on $C^+$   if and only if its Radon transform $R(\CF )\in D(\Bbb P^\vee )$ is micro-supported on $C^{\vee +}$.  \endproclaim 

\demo{Proof} Notice that  $C^{\vee +}= p^\vee_\circ p^\circ (C^+ )$ by 1.6.2, so for  $C^+\in \CC (\CF)$ one has  $ C^{\vee +}  \in \CC (R(\CF ))$ by 2.2. The converse  comes by interchanging $\Bbb P$ and $\Bbb P^\vee$ (since $C^+\in \CC (\CF)$ amounts to $C^+\in \CC(R^\vee R(\CF))$
 by  1.6.1(i) and 2.1(iii),(iv)).  \hfill$\square$
\enddemo

\subhead{\rm 3.4}\endsubhead  
Let $C$ be any strict  conical subset of $T^* \Bbb P$ and let $E\subset Q$ be the Legendre transform of $P(C)$.
Let $(h,f ): \Bbb P \leftarrow U \to Y$ be a test pair on $\Bbb P$. 
Let $p_U , p^{\vee }_U : Q_U :=U\times_{\Bbb P} Q \to U,\, \Bbb P^\vee$ be the projections; we have a test pair $(p^{\vee  }_U , fp_U ): \Bbb P^\vee \leftarrow Q_U \to Y$ on $\Bbb P^\vee$. Set $E_U := U\times_{\Bbb P} E \subset Q_U$.

\proclaim{Lemma}  $(h,f )$ is $C$-transversal if and only if $(p^{\vee }_U , fp_U \!)$ is $T^*\Bbb P^\vee$-transversal at  $E_U$, i.e., the map $ (p^{\vee }_U , fp_U \!) :  Q_U \to \Bbb P^\vee \times Y$ is smooth at $E_U$ (see Example (ii) in 1.2). 
\endproclaim

\demo{Proof } Pick any geometric point $u\in U$; set $x:=h(u)$, $y:= f(u)$. Then, by 1.6.2,
$C$-transversality of  $(h,f )$ at $u$ means that for every $e=(x,x^\vee )\in E_x$ the map $ d h_u +df_u : \lambda_{x^\vee x} \oplus   T^*_y Y \to T^*_u U$ is injective. The  latter assertion amounts to smoothness of $ (p^{\vee }_U , fp_U \!) :  Q_U \to \Bbb P^\vee \times Y$ at $(u,e)\in E_U$, i.e., to injectivity of the map
$dp_{U (u,e)}^\vee + d(fp_U )_{(u,e)} : T^*_{x^\vee} \Bbb P^\vee \oplus T^*_y Y \to T^*_{(u,e)}Q_U$.
To see this, notice that $T^*_{(u,e)}Q_U$ is the pushout of $T^*_u U \buildrel{dh}\over\leftarrow T^*_x \Bbb P \buildrel{dp}\over\to T^*_e Q$, hence it is the pushout of $T^*_u U \leftarrow \lambda_{x^\vee x} \to  T^*_{x^\vee} \Bbb P^\vee$ (see 1.6.2).
\hfill$\square$ \enddemo

 \demo{{\rm 3.5.} Proof of Theorem 3.2}  By 2.4, 1.6.1(i) and 2.1(iii),(iv) we can assume that $\CF =R^\vee (\CG )$, $\CG \in D(\Bbb P^\vee )$. Let $C \subset T^* \Bbb P$ be the cone over the Legendre transform of $E:= E_p (p^{\vee *}\CG )$.
   By 2.4 we need to show that  $C^+ = C^{\prime +}$ for every $C' \in \CCm (\CF )$.

Let $E' \subset Q$ be the Legendre transform   of $P(C')$. Then $
\Bbb P^\vee  \leftarrow Q\smallsetminus E' \to \Bbb P $ is a 
$C^{\prime \vee +}$-transversal test pair
 on $\Bbb P^\vee$ by 1.6.2. Now $\CG$ is micro-supported on $C^{\prime \vee +}$ by 3.3,  so, by  2.1(ii),  $p$ is universally locally acyclic relative to $p^{\vee *}\CG$ on $Q\smallsetminus E'$, i.e.,  $E \subset E'$ or $C^+ \subset C^{\prime +}$. By 2.4 it remains to prove that  $C^+ \in \CC (\CF )$. 
 
Let  $(h,f):\Bbb P \leftarrow U \to Y$ be a
  $C^+$-transversal test pair
 on $\Bbb P$. We need to check that it  is  $ \CF$-acyclic, i.e., $f$ is locally acyclic relative to $h^*\CF$.

 Consider the test pair $(p^{\vee  }_U , fp_U ): \Bbb P^\vee \leftarrow Q_U \to Y$ from 3.4. One has $h^*\CF =h^* p_* p^{\vee *}\CG = p_{U*}p^{\vee *}_U \CG$ by proper base change. Thus, since $p_U$ is proper, it is enough to show that  $ fp_U $ is locally acyclic   relative to $ p^{\vee *}_U \CG $ (see 3.13(i) below for the details), i.e., that $(p^{\vee  }_U , fp_U )  $ is $\CG$-acyclic.
We check this   separately (a) on $Q_U \smallsetminus E_U$, and (b) on a Zariski neighborhood of $E_U$:

(a) Since $p$ is universally locally acyclic relative to $p^{\vee *}\CG$ on $Q\smallsetminus E$ our $p_U$ is locally acyclic relative to $p_U^{\vee *}\CG$ on $Q_U \smallsetminus E_U $. Since $f$ is smooth (for $C^+$ contains the zero), $f p_U$ is locally acyclic relative to $p_U^{\vee *}\CG$ on $Q_U \smallsetminus E_U $    by \cite{D1} 2.14.  

(b) By 3.4,  $(p^{\vee  }_U , fp_U )$ is $T^* \Bbb P^\vee$-transversal at $E_U$ and so, by 1.2, on a neighborhood of $E_U$. Thus  $(p^{\vee  }_U , fp_U )  $ is $\CG$-acyclic by the lemma in 1.3. \hfill$\square$
\enddemo

\remark{Remark} Let $k'$ be any  finite extension of $k$. Consider the base change to $k'$ of the datum from 3.2. Clearly $E_{p_{k'} } (p_{k'}^{\vee *}R  (\CF_{k'} )  )=E_{p } (p^{\vee *}R  (\CF )  )_{k'}$. Thus $SS(\CF_{k'})=SS(\CF)_{k'}$. 
By 2.5(ii) and 2.3(i)  this holds for any smooth $X$ and $\CF\in D(X)$. \endremark

\subhead{\rm 3.6}\endsubhead It remains to prove the dimension estimate from   3.1. 
By 3.2 it can be rephrased as follows (we have interchanged $\Bbb P$ and $\Bbb P^\vee$ in 3.2 for notational convenience):

\proclaim{Theorem} For any  sheaf $\CG$ on $\Bbb P $  one has  $\dim E_{p^\vee} (p^* \CG )\le n-1$.  \endproclaim
 
We will  prove a slightly more general result with extra generic parameters (this generality is needed  for the induction argument). Let us formulate it.

\subhead{\rm 3.7}\endsubhead 
Suppose we have a relative version of the setting of 1.6, so there is an irreducible variety $S$ and a vector bundle $\CV$ on it; let $\Bbb P$ and $\Bbb P^\vee $ be the projectivizations of $\CV$ and the dual bundle  $\CV^\vee$, let $Q\subset \Bbb P\times_S \Bbb P^\vee$ be the incidence correspondence, and $p, p^\vee : Q \rightrightarrows \Bbb P,\,\Bbb P^\vee \!$ be the projections. 

\proclaim{Theorem } For every $\CG \in D(\Bbb P )$ there is a Zariski neighborhood $S^o$ of the 
generic point of $S$
such that $\dim E_{p^\vee} (p^* \CG )_{S^o} \le \dim \Bbb P -1$.
\endproclaim

The proof (inspired by Deligne's argument in  \cite{D1}) 
 takes the rest of the section. 
 
\remark{Remark} It is enough to prove the theorem with $S^o$ an {\it \'etale} neighborhood of the generic point (its image in $S$ is then the promised Zariski neighborhood). 
\endremark

\enspace

We fix a closed subset $D$ of $ \Bbb P$,  $D\neq \Bbb P$, such that $\CG$ is  locally constant on $\Bbb P\smallsetminus D$. Since $p^\vee$ is smooth, it is universally locally acyclic relative to $p^* \CG$   on
$Q\smallsetminus p^{-1}(D)$, i.e.,  $E_{p^\vee} (p^* \CG )\subset p^{-1}(D)$.

Denote the rank of $\CV$ by $n+1$, so $\dim \Bbb P=n+\dim S$. The proof  goes by induction by $n$. If 
$n=0, 1$  then $p$ and $p^\vee$ are isomorphisms, so the theorem   with $S^o=S$ follows directly from the above inclusion. From now on we assume that $n>1$.

\subhead{\rm 3.8}\endsubhead Let $T\subset \Bbb P$ be an $S$-family of lines (i.e., the projectivization of a rank 2 subbundle of $\CV$). Set $\Bbb P^\vee_{(T)}:=\{ x^\vee \in\Bbb P^\vee : Q_{x^\vee} \text{ is transversal to } T\}$; this is an open subset of $\Bbb P^\vee$ (its complement is an $S$-family of $n-2$-planes in $\Bbb P^\vee$). One has a smooth affine projection $ \Bbb P^\vee_{(T)} \to T$, $x^\vee \mapsto t=t(x^\vee ) := Q_{x^\vee}\cap T$. For an open subset $T^o$ of $T$ we denote by 
$\Bbb P^\vee_{(T^o )}$ the preimage of $T^o$ in $\Bbb P^\vee_{(T)}$ and by $Q_{(T^o )}$ the $ p^\vee $-preimage of   $\Bbb P^\vee_{(T^o )}$ in $Q$. We say that $T$ is {\it good} if $T\not\subset D$.

\proclaim{Proposition} For every good $T$ there is an open dense subset $T^o$ of $T$ such that $\dim E_{p^\vee} (p^* \CG )\cap Q_{(T^o )} \le \dim\Bbb P -1$.
\endproclaim

\demo{ Proposition  implies Theorem 3.7} To prove the theorem it is enough to find, after replacing $S$ by an \'etale neighborhood of its generic point, a collection $\{ U_\alpha \}$ of open subsets of $\Bbb P^\vee$ such that
$\dim E_{p^\vee} (p^* \CG )\cap p^{\vee -1}( U_\alpha )\le \dim\Bbb P -1$ and $\Bbb P^\vee \smallsetminus \cup U_\alpha$ is finite over $S$ (since due to the last property, one has then $\dim (Q \smallsetminus \cup \, p^{\vee -1}(U_\alpha ))\le \dim S + \dim Q/\Bbb P^\vee =  \dim\Bbb P -1$).  

Let us construct $U_\alpha$. We can assume that $\CV$ has a trivialization such that the coordinate axes sections $x_0 ,\ldots , x_n \in \Bbb P (S)$ do not lie in $D$.  
The $S$-lines  $T_{ij}$ that connect $x_i$ and $x_j$, $i\neq j$,  are good; let $T_{ij}^o$ be the open subsets as in the proposition. Shrinking $S$ we can assume that all $T_{ij}\smallsetminus T^o_{ij}$ are finite over $S$. The promised $U_\alpha$'s are $\Bbb P^\vee_{(T_{ij}^o )}$. (Here  the set 
$\Bbb P^\vee \smallsetminus \cup\Bbb P^\vee_{(T^o_{ij})}$ is finite over $S$ since any $x^\vee \in \Bbb P^\vee$ is uniquely determined by the datum of intersections $\{ Q_{x^\vee}\cap T_{ij} \} $).
 \hfill$\square$ \enddemo

The proof of the proposition is in 3.11; it is is based on the next two lemmas:

\subhead{\rm 3.9}\endsubhead  The first lemma is a general observation; its proof is in 3.12.  Suppose we have maps of varieties $ f: X\to Z$, $g: Y\to Z$,  $r: X\to Y$ such that
$f=gr$ and a sheaf $\CK \in D(X)$. 

\proclaim{Lemma} (i) If $r$ is proper, then $E_g (r_* \CK )\subset r(E_f (\CK))$.
\newline  (ii) If, in addition,
 $E_f (\CK )$ is finite over $Y$ then $E_g (r_* \CK )=r(E_f (\CK))$, hence one has $\dim E_g (r_* \CK )=\dim  E_f (\CK)$.
\endproclaim

\subhead{\rm 3.10}\endsubhead We return to the setting of 3.8; let $T$ be any  $S$-family of lines and let $T^o$ be any open subset of $T$. Consider a commutative diagram of incidence correspondences 
$$\spreadmatrixlines{1\jot}
\matrix  && \,\Bbb P &
\\
&q \nearrow\,\,&\!\pi  \uparrow &\nwarrow p 
\\
\Bbb P^\sim_{T^o}  &  \buildrel{p^\sim}\over\longleftarrow   &  Q^\sim_{(T^o )}
 &\buildrel{q^\sim}\over\lra  & Q_{(T^o )}
   \\  
r_T \downarrow && r \downarrow  &&  \downarrow   p^\vee
\\
 \bar{\Bbb P}_{T^o} &   \buildrel{\bar{p}_{(T^o )}}\over\longleftarrow &        \bar{Q}_{(T^o )} 
 &  \buildrel{\bar{p}_{(T^o )}^\vee}\over\lra  &  \Bbb P^\vee_{(T^o  )} 
  \\
 & \searrow& \downarrow& \swarrow
 \\
&& T^o
    \endmatrix
\tag 3.10.1$$
defined as follows. Below   $\bar{x}$ is a point of the Grassmannian $\text{Gr}=\text{Gr} (2,\CV )$ of lines in $\Bbb P$,  $L_{\bar{x}} \subset \Bbb P$ is the corresponding line, $t$ is a point of $T^o$, $x\in \Bbb P$,   $x^\vee \in \Bbb P^\vee_{(T^o )}$. 
A point of $\Bbb P^\sim_{T^o} $ is a triple $(t,x,\bar{x})$ with $t,x\in L_{\bar{x}} $. A point of $Q^\sim_{(T^o )}$ is a quadruple $(t,x, \bar{x},x^\vee )$ 
with $t,x \in L_{\bar{x}} \subset Q_{x^\vee}$. A point of $\bar{\Bbb P}_{T^o}$ is a pair $(t,\bar{x})$ with  $t\in   L_{\bar{x}}$. A point of $  \bar{Q}_{(T^o )} $ is a triple $(t, \bar{x},x^\vee )$ with $t\in L_{\bar{x}}\subset Q_{x^\vee}$. The arrows in (3.10.1) are the evident projections $q  ((t,x,\bar{x}))=x$, $p^\sim ((t, x,\bar{x},x^\vee ))= (t,x,\bar{x})$, etc.

\proclaim{Lemma} Suppose $T^o \subset T\smallsetminus D$. Then \newline (i)  $q^\sim$ yields an isomorphism of the closed sets $E_{p^\vee q^\sim}(\pi^* \CG ) \iso E_{p^\vee} (p^* \CG )\cap Q_{(T^o )}$. \newline (ii) $E_{p^\vee q^\sim}(\pi^* \CG )$ is finite over   $\bar{Q}_{(T^o )}$. 
\endproclaim

\demo{Proof of Lemma}  Since $p^\vee q^\sim$ and $p^\vee$ are smooth maps, one has  $E_{p^\vee q^\sim}(\pi^* \CG )\subset \pi^{-1}(D)$ and $
E_{p^\vee} (p^* \CG ) \subset p^{-1}(D) $. Now (i) follows since $q^\sim$ is an isomorphism over the open subset $x\neq t$ of $ Q_{(T^o )}$ which contains  $p^{-1}(D)\cap  Q_{(T^o )}$. To prove (ii) it is enough to check that $\pi^{-1}(D)$ is finite over  $\bar{Q}_{(T^o )}$. This comes since $r$ is a fibration by projective lines and  $\pi^{-1}(D)$ is a closed subset of $Q^\sim_{(T^o )}$ that does not contain any fiber of $r$ (indeed,   $\pi^{-1}(D)$   does not intersect the image of the section $(t, \bar{x},x^\vee )\mapsto (t,t, \bar{x},x^\vee )$ of $r$). \hfill$\square$
\enddemo

\subhead{\rm 3.11}\endsubhead {\it Proof of Proposition 3.8.} Our $T$ is good, so we have a nonempty $T^o$ as in the lemma in 3.10. Set 
$\bar{\CG }:= r_{T*}q^{  *}\CG \in D( \bar{\Bbb P}_{T^o}  )$.
One has
$\dim E_{p^\vee} (p^* \CG )\cap Q_{(T^o )}=\dim E_{p^\vee q^\sim}(\pi^* \CG )= \dim E_{\bar{p}^\vee_{(T^o )}}  (r_* \pi^* \CG  )= \dim E_{\bar{p}^\vee_{(T^o )}}  (\bar{p}_{(T^o )}^* \bar{\CG})$. Here the first equality comes from (i) of the lemma in 3.10; the second one comes from the lemma in 3.9 applied to $\CK =\pi^* \CG$, $f=p^\vee q^\sim$, $g=p^\vee_{(T^o )}$ (its conditions  hold due to (ii) of the lemma in 3.10); the last equality comes since $r_* \pi^* \CG = \bar{p}_{(T^o )}^* \bar{\CG}$ by the proper base change applied to $q^* \CG \in D(\Bbb P^\sim_{T^o})$ (the left commutative square in (3.10.1) is Cartesian and  $r_T$ is a proper map).

Now notice that 
 the bottom horizontal line in (3.10.1) is the standard incidence correspondence for the $T^o$-family of projective spaces  $\bar{\Bbb P}_{T^o}$  restricted to the open subset $\Bbb P^\vee_{(T^o )}$ of $\bar{\Bbb P}^\vee_{T^o}$ (the embedding $\Bbb P^\vee_{(T^o )} \hra \bar{\Bbb P}^\vee_{T^o}$ identifies $x^\vee$ with   $\bar{x}^\vee $ such that $\bar{Q}_{\bar{x}^\vee}= \bar{p}_{(T^o )}(\bar{p}_{(T^o )}^{\vee -1}(x^\vee ))$). Since $\dim \bar{\Bbb P}_{T^o} /T^o = n-1$,  the theorem in 3.7 is true for  $\bar{\Bbb P}_{T^o} /T^o$ and  $\bar{\CG }\in   D (\bar{\Bbb P}_{T^o} )$  by the induction hypothesis. Therefore, after shrinking $T^o$, we get $\dim E_{\bar{p}^\vee_{(T^o )}}  (\bar{p}_{(T^o )}^* \bar{\CG})\le \dim \bar{\Bbb P}_{T^o} -1 =\dim \Bbb P -1$, and we are done. \hfill$\square$

\subhead{\rm 3.12}\endsubhead 
 {\it Proof of Lemma 3.9.} Recall  the definition of  acyclicity (see \cite{D1} 2.12). Let $\CA (Z)$ be the collection of pairs $(z,h)$, $z$ is a geometric point of $Z$, $h$ is a geometric point of the Henselization $Z\tilde{_z}$ of $Z$ at $z$. For $f$, $\CK$ as above and $(z,h)\in \CA (Z)$ consider the maps of fibers $X_z \buildrel{i}\over\to X_{Z\tilde{_z}} \buildrel{\kappa}\over\leftarrow X_h$; we have the nearby cycles complex $\Psi_{f }(\CK )_{(z,h)}:=i^* R\kappa_* \CK_h$ (this  is a complex of \'etale sheaves on  $X_z$ that may not be constructible) and the evident map $\nu : \CK_z \to \Psi_{f }(\CK )_{(z,h)}$; the vanishing cycles complex 
$\Phi_{f }(\CK )_{(z,h)}$ is the cone of $\nu $. Then $f$ is  locally acyclic relative to  $\CK$  if  $\Phi_{f} (\CK )_{(z,h)}$ vanishes for every  $(z,h)\in\CA (Z)$,  and $f$ is universally locally acyclic if the same vanishing holds  for the base change $(f_{Z'}: X_{Z'}\to Z' ,\CK_{X_{Z'}})$ of $(f: X\to Z,\CK )$ by every  $Z'/ Z$. 

For $Z'/Z$ and  $(z',h')\in \CA (Z')$ let $S(Z'/Z,(z',h'), f, \CK)$ be the image  by $X_{z'}\to X$ of the support of the   complex $\Phi_{ f_{Z'}}(\CK_{X_{Z'}})_{ (z',h')}$. We see that $E_f (\CK )$ is the closure in $X$ of the union of subsets $S(Z'/Z,(z',h'), f, \CK)$ for all $Z'/Z$,  $(z',h')\in \CA (Z') $.

 For $g$, $r$ as above with $r$ proper  one has $\Phi_{g}(r_* \CK )_{(z,h)}=Rr_{z*} \Phi_{f} (\CK )_{(z,h)}$ by the  proper base change. The same is true after the base change by every $Z'/Z$. Thus $ S(Z'/Z,(z',h'), g, r_* \CK) \subset r(S(Z'/Z,(z',h'), f, \CK))$ and 3.9(i) follows since $r$ is proper.
 
  Suppose, in addition, that $r$ is  finite on $E_f (\CK )$. Each 
 $\Phi_f(\CK )_{(z,h)}$ is supported on $E_f (\CK )_z$. So, by the proper (or rather finite) base change, for every geometric point $y\in Y_z$ one has $ (\Phi_{g}(r_* \CK )_{(z,h)})_y = \oplus (\Phi_{f} (\CK )_{(z,h)})_x$ where $x$ runs the finite set   $E_f (\CK )_y$.  Hence the support of 
$  \Phi_{g} (r_* \CK )_{(z,h)}$ {\it equals} the $r_z$-image of the support of $ \Phi_{f}(\CK )_{(z,h)}$. The same is true after the base change by any $Z'/Z$. Therefore 
$ S(Z'/Z,(z',h'), g, r_* \CK) = r(S(Z'/Z,(z',h'), f, \CK))$ and 3.9(ii) follows since $r$ is proper. \hfill$\square$

\head 4. The singular support and the ramification divisor \endhead

  In this section we prove   the theorems from \S 1. The key one is the theorem  in 1.7. Here is an outline of its proof (minus the last assertion):
 
  We follow the notation from 1.7. It is enough to treat the case when $\CF$ is an irreducible perverse sheaf. Let $Y\subset \Bbb P$ be its support and let $C_\gamma$ be the irreducible components of $C:=SS(\CF )$. 
 By 3.1 one has   $\dim C \le n= \dim \Bbb P$. Then $D$ is the image of $P(i_\circ C)$ in $\tilde{\Bbb P}^\vee$. The geometry of the Veronese embedding shows that the images $D_\gamma$ of  $P(i_\circ C_\gamma )$ are distinct irreducible components of $D$ and the maps 
$P(i_\circ C_\gamma )\to D_\gamma$ are generically radicial
 (so $ \dim P(i_\circ C_\gamma )= N_d -n -1 + \dim C_\gamma  $ equals $    \dim D_\gamma$). One of the $C_\gamma$'s is the conormal bundle $T^*_Y \Bbb P$; it has dimension $n$, so the corresponding $D_\gamma$ is a hypersurface.  Let $\CG$ be the nonconstant irreducible constituent of $\tilde{R}(\CF )$, see 1.6.1(iii). Then  either (a) $SS(\CG )= (i_\circ C)^\vee$, or (b)  $SS(\CG )= (i_\circ C)^{\vee +}$ (see 3.3).  In case (a) the support of $\CG$ equals $D$; since it is irreducible, one has
 $C=T^*_Y \Bbb P$. In case (b) $\CG$ is the middle extension of an irreducible local system at the generic point of $\tilde{\Bbb P}$,  $D$ is its ramification locus, which is a divisor, and so each $C_\gamma$ has dimension $n$. 
 
 Below is the detailed story. We begin with preliminary lemmas (see 4.1--4.5).

 \subhead{\rm 4.1}\endsubhead    Let $X$ be a smooth variety.  For a closed subset $Y$ of $X$ we denote by $T^*_Y X$ a closed conical subset of $T^* X$ defined as follows. If $Y$ is generically smooth,\footnote{Which is always true  if the base field $k$ is perfect.} i.e., contains a dense open smooth (over $k$) subscheme $U_Y$, then $T^*_Y X$ is the  closure in $T^* X$ of the conormal bundle $T^*_{U_Y} X$ to $U_Y$. Otherwise one can find a finite  extension $k'$ (which we can choose to be purely inseparable) of $k$ such that $Y_{k'}\subset X_{k'}$ satisfies the above condition over $k'$, so we have $T^*_{Y_{k'}} X_{k'} \subset T^* X_{k'}=(T^*X)_{k'}$, and we define $T^*_Y X$ as the image of $T^*_{Y_{k'}}X_{k'}$ by the finite projection $T^* X_{k'}\to T^* X$. The above definition does not depend on the choice of $k'$. If $X$ is connected and $Y\neq \emptyset$, then $\dim T^*_Y X= \dim X$.
 
 \proclaim{Lemma}
 If $Y$ is the support of
 $\CF\in D(X)$,   then $T^*_Y X\subset SS(\CF )$.
  \endproclaim

\demo{ Proof} If $Y$ is generically smooth, then we can shrink $U_Y$ as above so that $\CF$ is locally constant on it, and the claim follows from 2.1(iii), 2.3(iii), and 2.5(ii). For a general $Y$ use Remark in 3.5 to reduce  to the generically smooth situation.
 \hfill$\square$
\enddemo

 \subhead{\rm 4.2}\endsubhead  Let $\pi : Q \to P$ be a  map of varieties; assume that $P$ is irreducible. For $a\ge 2$ we denote by $Q^{(a)}_P$ the complement to the union of all pairwise diagonals in the $a$-fold fiber product $Q_P^a$ of $Q$ over $P$.\footnote{So a geometric point of $Q^{(a  )}_P$ is a collection $(q_1,\ldots , q_a )$ of pairwise distinct geometric points of $Q$ such that $\pi(q_1 )=\ldots =\pi (q_a )$.} 
 We say that   $a$ (not necessary pairwise different) closed irreducible subsets $Z_1 ,\ldots ,Z_a$  of $Q$  {\it  intersect  well relative to} $\pi$ if the dimension of every irreducible component of $(Z_1 \times_P \ldots \times_P Z_a )\cap Q^{(a )}_P$ is equal to $  (\Sigma \dim Z_i ) - (a-1)\dim P = \dim P -\Sigma (\dim P -\dim Z_i )$. 

  Recall that a generically surjective map $ Z\to T$ with $Z$ irreducible is  {\it small} (in the Goresky--MacPherson sense) if $\dim (Z\times_T Z \smallsetminus \Delta (Z) )< \dim Z$; here $\Delta$ is the diagonal embedding. Notice that such a map is generically radicial over $T$.

 \proclaim{Lemma}  Suppose  $\pi : Q\to P$ is proper. \newline
 (i) If $Z\subset Q$ intersects well relative to $\pi$ with itself  and $\dim Z < \dim P$, then the map
 $ \pi|_Z : Z\to \pi (P)$  is small. \newline (ii) Suppose $Z_1 $, $Z_2$ intersect well relative to $\pi$, $Z_1$ is generically finite over $ \pi (Z_1 )$, and  $\dim Z_2 <\dim P$. Then
 $Z_1  \subset Z_2$ if (and only if) $ \pi (Z_1 ) \subset  \pi (Z_2 )$.   \endproclaim
 
 \demo{Proof}  (i) follows since $\dim (Z\times_P Z \smallsetminus \Delta (Z) )=\dim Z + (\dim Z -\dim P)$.
 Let us check (ii). Since $\dim (Z_1 \times_P Z_2 \smallsetminus \Delta (Z_1 \cap Z_2 ))= \dim Z_1 + (\dim Z_2 -\dim P)< \dim Z_1 =\dim   \pi (Z_1 )$, the fiber of $Z_1 \times_P Z_2 \smallsetminus \Delta (Z_1 \cap Z_2 )$ over the generic point of $  \pi (Z_1 ) $ is empty. Thus  $ \pi (Z_1 ) \subset  \pi (Z_2 )$
  implies that
 the fiber of $Z_1 \cap Z_2$ over the generic point of $ \pi (Z_1 )$ is nonempty. Hence $\dim (Z_1 \cap Z_2 )\ge \dim \pi (Z_1 )=\dim Z_1$, and so $Z_1 \cap Z_2 = Z_1$ since $Z_1$ is irreducible, i.e.,  $Z_1  \subset Z_2$.  \hfill$\square$
   \enddemo

 We say that a collection $\CS $ of closed irreducible subsets of $Q$  {\it intersects $m$-well relative to $\pi$},  $m\ge 2$, if every $a$ (not necessarily pairwise different) elements of $\CS$, where $m\ge a\ge 2$,   intersect well relative to $\pi$.  
 
 \subhead{\rm 4.3}\endsubhead We return to the situation of 1.7 and follow the notation there. So we consider  the Veronese embedding  $i= i_d : \Bbb P\hra \tilde{\Bbb P} $   of degree $d\ge 2$,  
 $\tilde{\Bbb P}^\vee =P(\tilde{V}^\vee)$, $\tilde{V}^\vee =\Gamma (\Bbb P,\CO (d))$. A closed conical subset $C$ of $T^* \Bbb P$ yields a strict 
closed conical  subset  $i_\circ C $ of $T^* \tilde{\Bbb P}$ (see 2.4), so we have $
 P(i_\circ C)\subset  P(T^* \tilde{\Bbb P})\iso \tilde{Q}$ and the 
projection $\tilde{p}^\vee_C : P(i_\circ C) \twoheadrightarrow D_C \subset \tilde{\Bbb P}^\vee$. If $C$ is irreducible, then so are $i_\circ C $,  $P(i_\circ C)$,  $D_C$.

 For a $\bar{k}$-point $x$   of $\Bbb P$, where $\bar{k}/k$ is any algebraically closed field, let $x^{(1)}$ be its first infinitesimal neighborhood   in $\Bbb P_{\bar{k}}$ and let  $J_{x }$ be the  vector space of sections of $\CO (d) $ on  $x^{(1)}$. Let $\rho_x :\tilde{V}^\vee_{\bar{k}}= \Gamma (\Bbb P_{\bar{k}}, \CO (d ))  \to J_x$ be the restriction map. 
 
 \enspace
 
 Assume that $m\ge 2$. Consider the following condition $(*)_{m,d}$:

 \enspace  

{\it For every collection of pairwise distinct  points $x_1, \ldots ,x_m$,  the product of the $\rho_{x_i}$ maps\footnote{In other words, $\rho_{(x_1 ,\ldots, x_m )}^{(m)}$ is the restriction map $\Gamma (\Bbb P, \CO (d ))_{\bar{k}}  \to \Gamma (\sqcup x_i^{(1)},\CO (d))$.}
$\rho^{(m)}_{(x_1 ,\ldots, x_m )}: \tilde{V}^\vee_{\bar{k}}  \to J_{x_1} \times \ldots \times J_{x_m} $ is surjective.} 

\enspace

Notice that for a given $m$  condition  $(*)_{m,d}$  is satisfied for large enough $d$.

\proclaim{Lemma}  $(*)_{m,d}$ implies that
 the collection of all subsets $\{ P(i_\circ C)\}$, $C$ is a closed conical subset of $T^* \Bbb P$, intersects $m$-well relative to $\tilde{p}^\vee :\tilde{Q}  \to  \tilde{\Bbb P}^\vee$.  \endproclaim

\demo{Proof} (i) Let us  describe $(P(i_\circ C_1 )\times_{\tilde{\Bbb P}^\vee}\ldots \times_{\tilde{\Bbb P}^\vee} P(i_\circ C_a ))\cap  \tilde{Q}^{(a )}_{\tilde{\Bbb P}^\vee}   $  explicitly.

Let $J$ be the vector bundle on $\Bbb P$ with fibers $J_x $; we have an evident exact sequence of locally free sheaves $0\to \Omega^1_{\Bbb P}(d)\to J \to \CO (d)\to 0$. Let $\Bbb P^{(a )}\subset \Bbb P^a$ be the complement to the union of  all pairwise diagonals in the product of $a$ copies of $\Bbb P$ and let $J^{a}_{\Bbb P^{(a)}}$ be the restriction of   the exterior product   of $a$ copies of $J$ to $\Bbb P^{(a)}$. Consider the map of vector bundles $\rho^{(a)}  :  \tilde{V}^\vee_{\Bbb P^{(a)}}\to J^{a}_{\Bbb P^{(a)}}$ whose fibers are $\rho^{(a)}_{(x_1 ,\ldots, x_a )}$. 

 Since $C_i$ are conical they yield closed conical subsets $C'_i$ in the $\CO_{\Bbb P}  (d)$-twisted cotangent bundle  $T^{\prime *} \Bbb P $ hence in $J$ (see above). The restriction $(C'_1 \times \ldots \times C'_a )_{\Bbb P^{(a)}}$ of the product of $C'_i$'s to  $ {\Bbb P^{(a)}}$ is a $\Bbb G_m^a$-invariant closed subvariety of $J^{a}_{\Bbb P^{(a)}}$. Thus
 $ \rho^{(a)-1}( (C'_1 \times \ldots \times C'_a )_{\Bbb P^{(a)}})$ is a closed conical subset of $\tilde{V}^\vee_{\Bbb P^{(a)}}$. Its projectivization equals 
 $(P(i_\circ C_1 )\times_{\tilde{\Bbb P}^\vee}\ldots \times_{\tilde{\Bbb P}^\vee} P(i_\circ C_a ))\cap  \tilde{Q}^{(a )}_{\tilde{\Bbb P}^\vee}$. 
 
(ii) Assume $(*)_{m,d}$ holds and $a\le m$. Then $ \rho^{(a)}$ is a surjective morphism of vector bundles. The dimension of its kernel is $ N_d +1 - a(n+1)$.
Below we view $ \rho^{(a)}$ as a smooth surjective map of smooth varieties whose fibers are affine planes of dimension $ N_d +1 - a(n+1)$.  
 
All irreducible components of  $(C'_1 \times \ldots \times C'_a )_{\Bbb P^{(a)}}$  have dimension $\Sigma \, \dim C_i$. Hence those of  $ \rho^{(a)-1}(
 (C'_1 \times \ldots \times C'_a )_{\Bbb P^{(a)}})$  have dimension $N_d +1 - a(n+1) +\Sigma  \,  \dim C_i = N_d +1 -  \Sigma \, (n+1 - \dim C_i )= N_d +1 -  \Sigma  \, (N_d +1 - \dim i_\circ C_i )$. Therefore, by (i), 
 all irreducible components of  $(P(i_\circ C_1 )\times_{\tilde{\Bbb P}^\vee}\ldots \times_{\tilde{\Bbb P}^\vee} P(i_\circ C_a ))\cap  \tilde{Q}^{(a )}_{\tilde{\Bbb P}^\vee}$ have dimension $\dim  \tilde{\Bbb P}^\vee - \Sigma \, ( \dim  \tilde{\Bbb P}^\vee   - \dim P( i_\circ C_i ))$. We are done.  \hfill$\square$
\enddemo
 
  \subhead{\rm 4.4}\endsubhead For our aims we need the case of $m=2$ in 4.3.
  
\proclaim{Lemma} Condition $(*)_{2,d}$ holds for every $d\ge 3$.
\endproclaim   

\demo{Proof} Induction by $n$. If $n=1$, then the assertion is evident. Suppose $n>1$. For distinct points $x_1$, $x_2$  of $\Bbb P=\Bbb P^n$, choose a hyperplane $H=\Bbb P^{n-1}$ passing through them. Let $J_{x_i}^{H}$ be the sections of $\CO (d)$ on the first infinitesimal neighborhood of $x_i$ in $H$. The exact sequence of sheaves $0\to \CO_{\Bbb P} (d-1)\to  \CO_{\Bbb P} (d)\to \CO_H (d)\to 0$ yields the short exact sequences  $0\to \Gamma (\Bbb P , \CO (d-1)) \to  \Gamma (\Bbb P , \CO (d))\to  \Gamma (H , \CO (d))\to 0$ and $0 \to \CO_{\Bbb P} (d-1)_{x_i} \to J_{x_i}\to J^H_{x_i}\to 0$. 
Since  the lemma for $H$ is known by the induction, we are done by the surjectivity of the restriction map $\Gamma (\Bbb P , \CO (d-1))\to \CO (d-1)_{x_1}\times \CO (d-1)_{x_2}$.  \hfill$\square$
\enddemo

  \subhead{\rm 4.5}\endsubhead    Assume that $d\ge 3$.  Combining  4.2--4.4 we get the next result:
 
 \proclaim{ Proposition}  (i)  If $C$ is any closed conical irreducible subset of $T^* \Bbb P$ of dimension  $ \le n$, then the map $\tilde{p}^\vee_C : P(i_\circ C)\to D_C$ is small, hence  generically radicial. \newline 
 (ii) Let $C$ be any closed conical subset of $T^* \Bbb P$ of dimension $\le n$. If $C_\gamma$ is an irreducible component of $C$, then $D_\gamma := D_{C_\gamma }$ is an irreducible component of $ 
D_C $, and  $C_\gamma \mapsto D_\gamma$ is 
a bijective correspondence  between the sets of irreducible components of $C$ and $D_C$. \hfill$\square$ \endproclaim
  
 The next observation will not be used in the sequel; the reader can skip it. For $(x,\nu )\in P(T^* \Bbb P )$ let  $P_{(x,\nu )} \subset \tilde{\Bbb P}^\vee$ be the linear projective  $(N-n)$-subspace formed by all degree $d$ hypersurfaces in $\Bbb P$ that pass through  $x\in\Bbb P$  and have the conormal direction $\nu\in P(T^*_x \Bbb P)$ at it.

\proclaim{Lemma}
 An irreducible hypersurface $D$ in $\tilde{\Bbb P}^\vee$  comes from some irreducible conical subset $C$ of dimension $n$ in $T^* \Bbb P$ if and only if $D$ is either the discriminant hypersurface (then $C$ is the zero section of $T^* \Bbb P$) or $D$ is spanned by all the subspaces   $P_{(x,\nu )}$ as above that lie in $D$.
\endproclaim 
 
\demo{Proof} If $D=D_C$, then $D$ is the union of the subspaces  $P_{(x,\nu )}$ for $(x,\nu )\in P(C)$ and of the projectivized conormals to $\Bbb P$ in $\tilde{\Bbb P}$ at points $x$ of the base of $C$. If $C$ is not the zero section of $T^* \Bbb P$, i.e., the base of $C$ is not $\Bbb P$, then the span of the latter subspaces has dimension $<N-1$, hence $D$ is spanned by   $P_{(x,\nu )}$.  

 Conversely, for an irreducible hypersurface  $D$, consider the closed subspace $K_D$ of $P(T^* \Bbb P )$
formed by those $(x,\nu )$ that $P_{(x,\nu )}$ lies in $D$.
If $D$ is spanned by  $P_{(x,\nu )}$, then $\dim K_D + (N-n)\ge \dim D=N-1$, i.e., $\dim K_D\ge n-1$. Pick any irreducible $P(C)\subset K_D$ of dimension $n-1$. Then the hypersurface $D_C$  lies in $D$, so $D_C =D$ by the irreducibility of $D$, and we are done.  \hfill$\square$ \enddemo

 \subhead{\rm 4.6}\endsubhead {\it Proof of Theorem 1.7 (except for its last assertion): case $d\ge 3$.} We follow the notation of 1.7.

 (i) {\it It is enough to prove the theorem when $\CF$ is an irreducible perverse sheaf}: Indeed, for any $\CF\in D(\Bbb P )$ let $\{ \CF_\alpha \}$ be the (perverse) Jordan--H\"older constituents of $\CF$, i.e., all irreducible perverse sheaves that occur in some ${}^p \CH^i \CF$. Suppose we know the theorem for all $\CF_\alpha$; let us prove it for $\CF$. Clearly $\{ i_* \CF_\alpha \}$ are the Jordan--H\"older constituents of $i_* \CF$. They are not locally constant; let $\CG_\alpha$ be the single   not locally constant Jordan--H\"older  constituent of $\tilde{R} (i_* \CF_\alpha )$, see 1.6.1(iii). Then  $\{ \CG_\alpha \}$ is the set of not locally constant Jordan--H\"older constituents of $\tilde{R} (i_* \CF )$  by 1.6.1(ii), so   $D_\CF =\cup_\alpha \, D_{\CF_\alpha }$.\footnote{We use the fact that on a smooth variety  
 every perverse subquotient of a locally constant perverse sheaf   is locally constant.} This implies all the properties of $D_\CF $ from the theorem. 
 
 It remains to check that $SS(\CF )=\cup_\alpha \, SS(\CF_\alpha )$. One has
 $SS(\CF )\subset \cup_\alpha SS(\CF_\alpha )$ by 2.1(iv). If the inclusion is strict,  then
  $SS(\CF )$ does not contain the generic point of an irreducible component $C_\gamma$ of some $SS(\CF_\alpha )$. Since $C_\gamma$ is uniquely determined by the corresponding divisor $D_\gamma$, we see that
the image of $  (i_\circ SS(\CF ))^\vee =  SS(i_* \CF )^\vee$ (see 2.5) in $\tilde{\Bbb P}^\vee$ does not contain the generic point of $D_\gamma$, i.e., $ \tilde{R}(i_* \CF )$ is locally constant there which contradicts $D_\CF =\cup_\alpha \,  D_{\CF_\alpha }$.

 (ii) For the rest of the proof we assume that our $\CF$ is an irreducible perverse sheaf on $\Bbb P$. Set $C:=SS (\CF )$; by 3.1 this is
 a closed conical subset of $T^* \Bbb P$ of dimension $\le n$. Let
 $C_\gamma$ be the irreducible components of $C$. By 2.5 one has $i_\circ C =SS(i_* \CF )$ so, by 2.4 and 3.3,  $SS(\tilde{R}(i_* \CF))$ equals either $(i_\circ C)^\vee$ or $(i_\circ C)^{\vee +}$. Thus $D=D_\CF  $ is the base of    $(i_\circ C )^\vee$ which is $D_C$. By   4.5  the map $C_\gamma \mapsto D_\gamma :=D_{C_\gamma }$ is a 1--1 correspondence between the irreducible components of $C$ and $D$, the projection $\tilde{p}^\vee_\gamma : P(i_\circ C_\gamma ) \to D_\gamma$ is small, and $C_\gamma$ is the only closed conical subset   of $T^* \Bbb P$ of dimension $\le n$ such that the $D_{C_\gamma}=D_\gamma$.
So to prove the theorem it remains to check that 
  $D$ is a divisor.

(iii) Since $\CF$ is irreducible its support $Y$ is also irreducible. Then 
$SS(\CF )$ contains $T^*_Y  \Bbb P$ (see 4.1).   Since $\dim T^*_Y  \Bbb P    =n$, our $T^*_Y  \Bbb P$ is an irreducible component of $C$.
 Let $\CG$ be the irreducible perverse sheaf subquotient of ${}^p \CH^{0} \tilde{R}(i_* \CF )$ which is not locally constant (see 1.6.1(iii)). Then
 $SS(\CG )$ equals either $(i_\circ C)^\vee$ or $(i_\circ C)^{\vee +}$ (see (ii)). So  the support of $\CG$ contains $Z:= D_{T^*_Y \Bbb P}$.  Since $\dim P (i_\circ T^*_Y \Bbb P )=N_d -1$,   4.5  implies that $Z$ is a hypersurface. Since $\CG$ is irreducible its support  is also irreducible, hence it equals either  $Z$ or  $\tilde{\Bbb P}^\vee$.

(iv) If the support of $\CG$ equals $Z$, then $D=Z$, so it is a divisor. Notice that since $D$ is irreducible then so is $C$ (see (ii)), i.e., $C=T^*_Y \Bbb P$.

If the support of $\CG$ is $\tilde{\Bbb P}^\vee$, then $\CG$ is the middle extension of a local system $\CG_\eta$ at the generic point of $\tilde{\Bbb P}^\vee$. Therefore $\tilde{\Bbb P}^\vee \smallsetminus D $ is  the maximal open subset to which $\CG_\eta$ extends as a local system, hence $D $ is a divisor (the ramification divisor of $\CG_\eta$). \hfill$\square$

\subhead{\rm 4.7}\endsubhead {\it Proof of Theorem 1.7 (except for its last assertion): case $d=2$.} In this section $i : \Bbb P\hra\tilde{\Bbb P}$ is the Veronese embedding of degree 2. We already know that all the components $C_\gamma$ of $C:=SS(\CF )$ have dimension $n$. Since $D_\CF $ is the base of $(i_\circ C)^\vee $  by 2.5 and 3.3, the theorem follows from 4.2 and the next result:

\proclaim{Proposition} Let $C_1$, $C_2$ be any irreducible conical subsets of $T^* \Bbb P$ of dimension $n$. Then the subsets $P(i_\circ C_1 )$, $P(i_\circ C_2 )$ of $P(T^* \tilde{\Bbb P})=\tilde{Q}$ intersect well relative to $\tilde{p}^\vee$ (see 4.2) unless $C_i $ coincide and equal to the conormal bundle 
$T^*_{\Bbb P'} \Bbb P$ to some linearly embedded projective subspace $\Bbb P' $ of $ \Bbb P$ of dimension $  \ge 1$. The projection of $P(i_\circ T^*_{\Bbb P'}\Bbb P )$ to its image in $\tilde{\Bbb P}^\vee$ is   generically radicial. 
\endproclaim

To prove the proposition, we need the next lemma. We  use the notation from 4.3 for $d=2$. Let  $x_1$, $x_2$ be two distinct $\bar{k}$-points of $\Bbb P$.

\proclaim{Lemma} The image   of    
$\rho^{(2)}_{(x_1 ,x_2 )}:  \tilde{V}^\vee_{\bar{k}} \to J_{x_1}\times J_{x_2} $ has codimension 1. Its intersection with $T_{x_1}^{\prime *}\Bbb P\times T_{x_2}^{\prime *} \Bbb P$ is the preimage of the diagonal
$\bar{k} \subset \bar{k}\times \bar{k}$ by the restriction map 
$T_{x_1}^{\prime *}\Bbb P\times T_{x_2}^{\prime *} \Bbb P \twoheadrightarrow 
T_{x_1}^{\prime *} L \times T_{x_2}^{\prime *} L =  \bar{k}\times \bar{k}$ where $L= L(x_1 ,x_2 )$ is the projective  line that connects $x_1$ and $x_2$. Here the latter identification comes since $\Omega^1_{L}(2)=\CO_{\Bbb P^1}$. 
\endproclaim

\demo{Proof} The case $n=1$ is checked directly; then use the induction argument from the proof in  4.4. \hfill$\square$
\enddemo

\demo{Proof of Proposition} Let $Y_i \subset \Bbb P$ be the base of $C_i$. 
Suppose that $P(i_\circ C_1 )$ and $P(i_\circ C_2 )$  do {\it not} intersect well relative to $\tilde{p}^\vee$.\footnote{Notice that this evidently excludes the situation when $\dim Y_1=\dim Y_2 =0 $.}
A simple modification of  the proof  in 4.3 combined with the codimension 1 assertion from the lemma  shows that this happens
 if and only if, for all pairs of distinct  points $x_1 \in Y_1$, $x_2 \in Y_2$,  
   the product of the fibers $C'_{1 x_1}\times C'_{2 x_2} \subset T_{x_1}^{\prime *}\Bbb P \times T_{x_2}^{\prime *}\Bbb P$ is  contained in the image of $\rho^{(2)}_{(x_1 ,x_2 )}$. Since this product is $\Bbb G_m \times \Bbb G_m$-invariant, by the lemma,  this amounts to  the fact that the images of $C_{i x_i}$ in $T_{x_i}^{ *} L(x_1 ,x_2 )$ are both equal to 0. Which means that $C_{1 x_1}\subset T^*_{x_1}\Bbb P$ lies in the orthogonal complement to the vector subspace  $F_{1x_1}$ of $T_{x_1}\Bbb P$ generated by the tangents at $x_1$ to   lines $L(x_1 ,x_2 )$ for all  $x_2 \in Y_2$, plus the same condition   with indices 1 and 2 interchanged. Notice that: \newline
(a) One has $\dim C_{i x_i}   \le \dim F_{i x_i}^\perp  $ and the equality means that $C_{i x_i} = F_{i x_i}^\perp$. \newline (b) One has
$\dim Y_2 \le \dim F_{1x_1} $ and the equality means that $Y_2$ coincides with the linear projective subspace of $\Bbb P$ spanned by  the  lines  $L(x_1 ,x_2 )$ for all $x_2 \in Y_2$, i.e., $Y_2$ is a linear projective subspace of $\Bbb P$ that contains $x_1$. The same  holds with indices 1 and 2 interchanged.

If $x_i$ are general enough, then $\dim C_{i x_i} + \dim Y_i=\dim C_i =n$.
Since $\dim F_{i x_i}  +\dim F_{i x_i}^\perp  =n$, we can rewrite the inequality in (a) as $ \dim F_{i x_i}\le \dim  Y_i $. Combining it with (b) we see that the inequalities in (a) and (b) are equalities. Therefore $Y_i$ coincide and are equal to a linear projective subspace $\Bbb P'$ of $\Bbb P$, $\dim \Bbb P' >0$, and $C_i$ are both equal to  the conormal bundle $T^*_{\Bbb P'}\Bbb P$.

It remains to check the last assertion of the proposition. A point of $P(i_\circ T^*_{\Bbb P'}\Bbb P)$ is a pair $(x, \tilde{x}^\vee )$, $x\in \Bbb P'$, $\tilde{x}^\vee \in \tilde{\Bbb P}^\vee$, such that either $Q'_{\tilde{x}^\vee}:= 
\tilde{Q}_{\tilde{x}^\vee}\cap \Bbb P'$ equals $\Bbb P'$ or this is a quadric that  contains $x$ as a singular point. The projection to $\tilde{\Bbb P}^\vee$ is $(x, \tilde{x}^\vee )\mapsto \tilde{x}^\vee$.
We are done since  a generic singular quadric has single singular point.
\hfill$\square$
\enddemo

\remark{Remark}  For $\Bbb P'$ as above the dimension of
$P(i_\circ T^*_{\Bbb P'}\Bbb P)\times_{\tilde{\Bbb P}^\vee}P(i_\circ T^*_{\Bbb P'}\Bbb P)\smallsetminus \Delta (P(i_\circ T^*_{\Bbb P'}\Bbb P))$ equals $ \dim P(i_\circ T^*_{\Bbb P'}\Bbb P)$, so the projection of $P(i_\circ T^*_{\Bbb P'}\Bbb P )$ to its image in $\tilde{\Bbb P}^\vee$ is semi-small. 
 \endremark

 \subhead{\rm 4.8}\endsubhead {\it A   reminder about pencils.}
 We need a variant of  the classical Lefschetz pencils story from  \cite{Katz}. We are in the setting of 1.7 and follow the notation there. Let $C$ be any irreducible closed conical subset of $T^* \Bbb P$ of dimension $n$; set $D:=D_C$.
 
 (i) For a line $L$   in $\tilde{\Bbb P}^\vee$ set $\tilde{Q}_{L}:=\tilde{Q} \times_{\tilde{\Bbb P}} L \subset \tilde{Q}$; let $\pi_L :  
 \tilde{Q}_{L}\to L$ be the projection. 
Let $\tilde{U}_L \subset \tilde{\Bbb P}$ be the complement to the axis $L^\perp$ of $L$. The map $\tilde{p}|_{\tilde{Q}_{L}}:  \tilde{Q}_{L} \to \tilde{\Bbb P}$ is an isomorphism over $\tilde{U}_L $ so we have the inverse open embedding $j_L : \tilde{U}_L \hra \tilde{Q}_{L} $ and hence $\tilde{f}_L := \pi_L j_L : \tilde{U}_L \to L$. The map $\tilde{f}_L $ is smooth;  the
embedding $\tilde{U}_L \hra  \tilde{Q}_{L}\hra \tilde{Q}\iso P(T^* \tilde{\Bbb P})$ sends $\tilde{x}\in \tilde{U}_L $ to the conormal to the fiber of $\tilde{f}_L $ passing through $\tilde{x}$.

(ii) Set $U_L := i^{-1}(\tilde{U}_L )\subset \Bbb P$, $f_L := \tilde{f}_L  i$. Let $x$ be point in $U_L$. Pick a linear coordinate $t$ on $L$ regular at $f_L (x)$. The composition  $\tilde{f}_{(L,t)} :=t(\tilde{f}_L )$ is a function on a neighborhood $\tilde{U}$ of $i(x)$ in $\tilde{U}_L$, so $f_{(L,t)} :=\tilde{f}_{(L,t)} i$ is a function on $U :=i^{-1}(\tilde{U} )$. One has    $P(i_\circ C)_{\tilde{U} }\cap  \tilde{Q}_{L} \buildrel{\sim}\over{\leftarrow} d\tilde{f}_{(L,t)} (\tilde{U} )\cap i_\circ C \iso df_{(L,t)} (U  )\cap C$.

(iii) Pick a hyperplane $H\subset \Bbb P$ that does not contain $x$. We refer to   functions on  $\Bbb P \smallsetminus H$ as polynomials. Writing $\CO_{\Bbb P}(d) \iso \CO_{\Bbb P}(dH)$ we identify  $\Gamma (\Bbb P,\CO_{\Bbb P}(d))$ with the vector space of polynomials of degree $\le d$, so $\tilde{\Bbb P}^\vee$ is its projectivization. One has
$f_{(L,t)}= q_1 /q_2$ where $q_1 ,q_2$ are  nonzero polynomials of degree $\le d$ such that  $(q_1 ), (q_2 )\in L\subset  \tilde{\Bbb P}^\vee$, $t((q_1 ))=0$, $t((q_2 )) =\infty$. Conversely, if $q_1 ,q_2$ are any non-proportional nonzero polynomials of degree $\le d$ and $q_2 (x)\neq 0$,   then $q_1 /q_2 \in \CO_{\Bbb P x}$ equals $f_{(L,t)} $ for some $(L,t)$ as in (ii).

(iv) Consider the closed conical subset  $(i_\circ C)^\vee$ of $T^* \tilde{\Bbb P}^\vee$ of dimension $N$ with base $D$ (see 3.3). Let $\tilde{x}^\vee$ be a closed point of  $L\cap D$. Then
the embedding $L\hra \tilde{\Bbb P}^\vee$ is $(i_\circ C)^\vee$-transversal at $\tilde{x}^\vee$ if and only if the fiber $P(i_\circ C)_{\tilde{x}^\vee}$ lies in $U_L$. 
Since $ U_L$ is affine this implies that $P(i_\circ C)\to D$ is finite hence radicial at $\tilde{x}^\vee$.

 \subhead{\rm 4.9}\endsubhead {\it Proof of Theorem 1.5.} 
  By 2.5 and 2.3(i)  the case of general $X$ is reduced  to the case $X=\Bbb A^n$ then to $X=\Bbb P$. So we are in the setting of 1.7, and  we use the notation from there.
 
(i) Since $SS^w  (\CF )\subset SS (\CF )$, to prove the equality it is enough to find for every component $C_\gamma$ of $C=SS(\CF )$ and its open dense $\Bbb G_m$-invariant subset $W$ a datum $(x,U,f)$ where $ x$ is a closed point of $\Bbb P$, $U$ is a Zariski neighborhood of $x $ in $\Bbb P$,\footnote{Here we assume that $k$ is infinite; otherwise $ U\subset  \Bbb P_{k'}$ for a finite extension $k'$ of $k$, see 1.5.} and $f$ is a function on $U$  such that   $df (x)\in W_x$  and  $\phi_{f} ( \CF )_x \neq 0$.

(ii) Recall that  $D_\gamma :=D_{C_\gamma}$  is a component of the divisor $D=D_\CF$ and the map $\tilde{p}^\vee_\gamma : P(i_\circ C_\gamma )\to D_\gamma$ is generically radicial. 
Our $W$ yields an open dense subset $P(i_\circ C_\gamma )_W $ in $P(i_\circ C_\gamma )$ defined as the image of $W$
by the correspondence $T^*\Bbb P \buildrel{di}\over\longleftarrow
T^* \tilde{\Bbb P}|_{\Bbb P}\smallsetminus \{ \text{zero section} \} \to P(T^* \tilde{\Bbb P}|_{\Bbb P})\subset P(T^* \tilde{\Bbb P})$.  Choose an  open dense subset $D_\gamma^o$   of $D_\gamma$ such that
$P(i_\circ C_\gamma )^o :=\tilde{p}^{\vee -1}_\gamma (D_\gamma^o )\to D^o_\gamma $ is radicial, $D_\gamma^o$ does not intersect other components of $D$, the sheaf $\tilde{R} (i_* \CF  )|_{D^o_\gamma}$ is locally constant, and 
$P(i_\circ C_\gamma )^o \subset  P(i_\circ C_\gamma )_W$.

By Bertini the lines $L$ in $\tilde{\Bbb P}^\vee$ that   properly intersect   $D_\gamma$ at points of $ D^o_\gamma$ and such that the embedding $L\hra \tilde{\Bbb P}^\vee$ is $(i_\circ C_\gamma )^\vee$-transversal form a dense open subset of the Grassmannian of lines. So  one can pick one such
 $L$  defined over $k$ if $k$ is infinite or over a finite extension $k'$ of $k$ if $k$ is finite. Pick a point $\tilde{x}^\vee$ in $L\cap D_\gamma =L\cap D_\gamma^o$ and let $(i(x), \tilde{x}^\vee )$ be its preimage in $P(i_\circ C )^o$. Due to the
 $(i_\circ C_\gamma )^\vee$-transversality, $i(x)$ does not lie in the axis of $L$, i.e., $x\in U_L$ (see 4.8(iv)). Pick $t$ and a neighborhood $U\subset U_L$ of $x$ as in 4.8(ii), and set $f:=f_{(L,t)}$. Then $(x,U,f)$ is the promised datum from (i). Indeed, $df(x)\in W_x$ since $P(i_\circ C_\gamma )^o \subset  P(i_\circ C_\gamma )_W$ (see 4.8(ii)). It remains to check that   $\phi_{f} ( \CF )_x \neq 0$.

(iii) Let $L^o  $ be the complement in $L$ to the finite set  $(D\cap L)\smallsetminus \{ x^\vee\}$; set $\tilde{Q}_{L^o}:=\tilde{Q}_{L}\times_L L^o =  \tilde{Q} \times_{\tilde{\Bbb P}} L^o$.  
 Since $\tilde{p}^\vee :   \tilde{Q} \to \tilde{\Bbb P}^\vee$ is universally locally acyclic relative to  $\tilde{p}^* i_* \CF$ on $\tilde{Q}\smallsetminus P(i_\circ C )$, we see that the projection $\pi_{L^o}: \tilde{Q}_{L^o}\to L^o$
 is locally acyclic relative to $(\tilde{p}^* i_* \CF )|_{\tilde{Q}_{L^o}}$  on the complement to $\{ (i(x),\tilde{x}^\vee )\}= \tilde{Q}_{L^o} \cap P(i_\circ C ) $. Since  $\tilde{R} (i_* \CF )|_{L^o}=\pi_{L^o *}((\tilde{p}^* i_* \CF )|_{\tilde{Q}_{L^o}})$ by the proper base change and $\pi_{L^o}$ is proper we see that $\phi_{f} (\CF )_x  $ is equal to the vanishing cycles of $\tilde{R}(i_* \CF )|_L$ at $x^\vee$. The latter is nonzero for $\tilde{R}(i_* \CF )|_L$  is   {\it not} locally constant at $x^\vee$ (which follows since
 $ \tilde{R}(i_* \CF )|_{D_\gamma^o}$ is locally constant but $  \tilde{R}(i_* \CF )$ is {\it not} locally constant at $D_\gamma^o$). We are done. \hfill$\square$

\subhead{\rm 4.10}\endsubhead {\it Proof of the theorems 1.3(ii) and 1.4.} Assertions 1.3(ii) and 1.4(iii) for general $X$ reduce, using 2.5(ii) and 2.3(i) (or 2.5(i) and 1.5), to the case of 
 $X=\Bbb P$.  Here they follow from the description of $SS(\CF )$ of the theorem in 1.7.  To check 1.4(ii) we can replace, by the theorem in 1.5 (see 4.9), $SS$ by $SS^w$; now the claim follows from perverse t-exactness of functors $\phi_f$. Finally let us prove
1.4(i). One has $SS(\pi^* \CF )\subset \pi^\circ SS (\CF )$ by 2.2(i). As above, we prove the inverse inclusion with $SS$ replaced by $SS^w$. Now (see 1.5)  $SS^w (\CF )$ is the closure of the set of points $(x,df(x))\in T^* X$ where $x\in X$ and $f$ is a function on a neighborhood of $x$ which is not locally acyclic relative to $\CF$ at $x$ (with changes as in loc.cit.~for finite $k$). If $(x,df(x))$  is as above, $z\in \pi^{-1}(x) $,   then  $f\pi$ is not locally acyclic relative to $\pi^* \CF$ at $z$ (since $\pi$ is smooth) hence $(z, d(f\pi )(z)) \in SS^w (\pi^*\CF )$. We are done for $\pi^\circ   SS(\CF )$ is the closure of the set of such points $(z, d(f\pi )(z))$.\footnote{Use the fact that for any subset   $A\subset T^* X$ one has
 $ \pi^\circ \overline{A}=\overline {\pi^\circ A} $ where $\bar{\,\,}$ means the closure  (which follows since $\pi$ is open and 
 $d\pi : T^* X\times_X Z \to T^* Z$ is a closed embedding).}   \hfill$\square$

\subhead{\rm 4.11}\endsubhead {\it A linear algebra lemma.} The aim of the rest of the article is to prove  the last assertion of  Theorem 1.7. We assume that the base field $k$ is algebraically closed.

\enspace 

Let $W$ be a $k$-vector space of dimension $2n$, $\omega$ a symplectic (i.e., alternate non-degenerate) form on $W$. If $p$ (the characteristic of $k$) equals  $2$ then let $\kappa$ be a quadratic form on $W$ whose polarization equals $\omega$. 

A vector subspace $\CL$  of $W$ is said to be {\it $\omega$-isotropic}, resp.~{\it $\kappa$-isotropic}, if $\omega$, resp.~$\kappa$, vanishes on it. Denote by $\text{Gr}^\omega =\text{Gr}^\omega (W)$, resp.~$\text{Gr}^\kappa =\text{Gr}^\kappa (W)$, the Grassmanians of $\omega$-isotropic, resp.~$\kappa$-isotropic subspaces of $W$ of dimension $n$. The Grassmannian $\text{Gr}^\omega$ is irreducible, and $\text{Gr}^\kappa \subset \text{Gr}^\omega$ has 2 connected components: the subspaces $\CL^\kappa_1 , \CL^\kappa_2 \in \text{Gr}^\kappa$ lie in the same component if and only if   $\dim \CL^\kappa_1 /(\CL^\kappa_1 \cap \CL^\kappa_2 )$ is even.

\proclaim{Lemma} Let $V$ be any vector subspace of $W$ of dimension $n$. \newline (i) There is $\CL^\omega \in \text{\rm Gr}^\omega$ complementary to $V$.  \newline (ii) For $p=2$ there is $\CL^\kappa \in \text{\rm Gr}^\kappa$ complementary to $V$. \newline (iii) If  $V$ is not $\kappa$-isotropic, then  $\CL^\kappa$'s as in (ii) occur in both components of $\text{\rm Gr}^\kappa$.
\endproclaim

\demo{Proof} (i) Let us construct $\CL^\omega$. Set $V_0 :=V\cap V^\perp$ (here ${}^\perp$ is the orthogonal complement for $\omega$) and pick $V_1 \subset V$ complementary to $V_0$. Pick an $\omega$-isotropic $V'_0 \subset V_1^\perp$ such that the $\omega $-pairing between $V_0$ and $V'_0$ is nondegenerate. Set $V_2 := (V\oplus V'_0 )^\perp$. Pick an isomorphism $g:V_1 \iso V_2$ such that $g^* \omega|_{V_2}=-\omega|_{V_1}$ (it exists since $\dim V_1 =\dim V_2$ and $\omega$ is nondegenerate on $V_2$ and $V_1$); let $G\subset V_1 \oplus V_2$ be the graph of $g$. Set $\CL^\omega := V'_0 \oplus G$.

(ii) Pick an $\CL^\omega$ as in (i). If it is $\kappa$-isotropic we are done. If not then, since $\omega$ identifies $V$ with the dual to $\CL^\omega$, there is a nonzero $v \in V$ such that for  $l\in \CL^\omega$ one has $\kappa(l ) = \omega (v,l )$. 
Set $\CL^\omega_0 := v^\perp \cap \CL^\omega$, pick any $l\in \CL^\omega \smallsetminus  \CL^\omega_0$, and let $P$ be the plane generated by $v$ and $l$. Since $\omega|_P$ is nondegenerate,  $\kappa$ vanishes on exactly two lines in $P$; let $\ell$ be one of these lines that does not contain $v$. Set $\CL^\kappa := \ell \oplus \CL^\omega_0$.

(iii) We construct $\CL^\kappa_1$, $\CL^\kappa_2$ as in (ii) such that $\dim ( \CL^\kappa_1 /\CL^\kappa_1 \cap \CL^\kappa_2 )=1$ so they belong to different components of $\text{Gr}^\kappa$. Pick an $\CL^\kappa$ as in (ii).
Pick $v\in  V$ such that $\kappa (v)\neq 0$.  Set $\CL^\kappa_0 := v^\perp \cap \CL^\kappa$, pick any $l\in  \CL^\kappa \smallsetminus \CL^\kappa_0$,
and let $P$ be the plane generated by $v$ and $l$. Since $\omega|_P$ is nondegenerate,  $\kappa$ vanishes on exactly two lines $\ell_1$ and $\ell_2$ in $P$. Set $\CL^\kappa_i := \ell_i \oplus \CL^\kappa_0$.  \hfill$\square$
\enddemo 

\remark{Remark}  $\CL^\omega$ and $\CL^\kappa$ as in the lemma form Zariski open subsets of $\text{Gr}^\omega$ and $\text{Gr}^\kappa$. By the lemma, they are dense unless $p=2$ and $V$ is $\kappa$-isotropic and we deal with $\text{Gr}^\kappa$: here the closure of our open subset is one of the two components of $\text{Gr}^\kappa$.
\endremark

\subhead{\rm 4.12}\endsubhead  For a smooth  variety $X$, $\dim X=n$, consider  the tangent bundle $T (T^* X)$. Let $T^v (T^* X ):= T (T^* X /X)$ be its vertical subbundle, so for  $(x,\nu )\in T^* X$ one has an evident canonical identification $T^*_x X \iso  T^v_{(x,\nu)} (T^* X )$, $\mu \mapsto \mu^v$. A function $f$ on an open subset $U$ of $X$ yields a section $df : U\to T^* X$; we denote by $d^{(2)}f : TU \to  T (T^* X)$ its differential. So for   $(x,\nu )\in T^* X$, every $f$ with $df (x)=\nu$ yields a direct sum decomposition $\alpha_{f,x}: T^*_x X \oplus T_x X \iso T_{(x,\nu)} (T^* X )$, $\alpha_{f,x} (\mu ,\tau ):= \mu^v + d^{(2)}f (\tau )$.

Recall that $T^* X$ is naturally a symplectic manifold, so $T (T^* X)$ carries a canonical symplectic form $\omega$. It vanishes on $T^v (T^* X )$ and on the image of $d^{(2)}f $ for every $f$. Explicitly,  $\alpha_{f,x}$ identifies $\omega_{(x,\nu )}$ with the standard symplectic form on $T^*_x X \oplus T_x X$.

If $p=2$, then $T (T^* X)$ carries a  canonical quadratic form $\kappa$ such that the polarization of $\kappa$ equals $\omega$ and $\kappa$ vanishes on the image of $d^{(2)}f$ for every   $(U, f)$ as above.  Our $\kappa$ is uniquely defined by these properties; it vanishes on $T^v (T^* X)$. Explicitly, one has  $\kappa  (\alpha_{f,x} (\mu ,\tau ))= \mu(\tau )$. 

\proclaim{Lemma} Let $\CL$ be a vector subspace of $T_{(x,\nu )} (T^* X )$
 complementary to $T^v_{(x,\nu )} (T^* X )$. Then $\CL$
   can be realized as the $d^{(2)}f$-image of $T_x X$ for some $f\in \CO_{X x}$,  $df(x)=\nu$, if (and only if)  $\CL$ is $\omega$-isotropic for $p\neq 2$  and $\CL$ is $\kappa$-isotropic for  $p=2$.
\endproclaim

\demo{Proof}  We pick local coordinates with $x=0$ and  identify functions on $T_x X$ with polynomial functions of coordinates.  Our  $f$ is the sum of a linear term equal to $\nu$ and a quadratic function $q$. Then $d^{(2)}f$ is  the graph of the polarization of $q$.

Our $\CL$ is the graph of a bilinear function $B$. The condition of the lemma means that $B$ is alternative, and so $B$ is the polarization of some $q$. We are done.         \hfill$\square$
\enddemo

 Let $C$ be an irreducible closed conical subset of $T^* X$ of dimension $n$ and let  $(x,\nu )$ be a closed point of $ C_{\text{reg}}$.
 For $f\in \CO_{Xx}$ such that $(x,\nu )$ is an isolated intersection point of    $C$ and $df(U)$,\footnote{Here $U$ is a neighborhood of $x$ where $f$ is defined.} let $\langle C, df  \rangle_{(x,\nu )}$ be the intersection index of these subvarieties at $(x,\nu )$ in $T^* X$. Thus $\langle C, df  \rangle_{(x,\nu)}=1$ means that the intersection is transversal.
 
 \proclaim{Proposition}\!\!\!\footnote{Inspired by a discussion with Deligne.} (i)  One can find $f\in\CO_{Xx}$ with $\langle C, df  \rangle_{(x,\nu )}=1$ if and only if the next condition $(*)$ is not satisfied: \newline  $(*)$ $p=2$, $T_{(x,\nu )}C$ is $\kappa$-isotropic, and the rank of the map $T_{(x,\nu )}C \to T_x X$ is odd.
 \newline
(ii)  In the situation of $(*)$  one can find $f\in\CO_{Xx}$ with $\langle C, df  \rangle_{(x,\nu )}=2$.    \endproclaim

\demo{Proof} 
(i) By the lemma one can find $f$ as in (i) if and only if there is
an $\omega$-isotropic, or $\kappa$-isotropic if $p=2$, subspace $\CL$ of  $T_{(x,\nu )}(T^* X)$  complementary to both $T_{(x,\nu )} C$ and $T^v_{(x,\nu )}(T^* X)$. Let us try to find such an $\CL$.
 
By Remark in 4.11 applied to $W= T_{(x,\nu )}(T^* X)$, our $\CL$   always exists unless $p=2$ and $T_{(x,\nu )} C$ is $\kappa$-isotropic. In the latter situation  it exists if and only if $T_{(x,\nu )} C$ and $T^v_{(x,\nu )}(T^* X)$ lie in the same component of $\text{Gr}^\kappa$, i.e., $  T_{(x,\nu )} C/(T_{(x,\nu )} C \cap  T^v_{(x,\nu )}(T^* X))$  has even dimension. We are done.

(ii) Suppose $(*)$ holds so $T_{(x,\nu )} C$ and $T^v_{(x,\nu )}(T^* X)$ lie in   different components of $\text{Gr}^\kappa$.  
Then $\CL^\kappa$'s with  $\dim (\CL^\kappa \cap T^v_{(x,\nu )}(T^* X))=0$ and with 
$\dim (\CL^\kappa \cap T_{(x,\nu )} C )=1$
form dense open subsets  of the same component of $\text{Gr}^\kappa$. So there is  $\CL^\kappa$ that satisfies  both conditions; pick it. Set $\ell :=
\CL^\kappa \cap T_{(x,\nu )} C$.

Let $g_1,\ldots ,g_n \in \CO_{T^* X \, (x,\nu )}$ be local equations of $C$ near $(x,\nu )$ so $dg_i ((x,\nu))$ form a base of the conormal to $C$ at $(x,\nu )$. We can assume that $dg_n  ((x,\nu))$ vanishes on $\CL^\kappa$. Then the restriction of $\{ dg_i ((x,\nu))\}_{i\le n-1}$ to $\CL^\kappa$ is a base of the dual to $\CL^\kappa / \ell $. Since $dg_n  ((x,\nu))$ vanishes on both Lagrangian subspaces
 $\CL^\kappa $ and $T_{(x,\nu )} C$ one has
$dg_n  ((x,\nu))=\omega (\tilde{\tau} ,\cdot )$ for some generator $\tilde{\tau}$ of the line $\ell$; let $\tau$ be its image in $T_x X$.

Consider the set $\CS$ of $f\in \CO_{Xx}$ with $df  (x)=\nu$ and $d^{(2 )}f  (T_x X)= \CL^\kappa$; by the lemma it is not empty. For $f\in \CS$ set $r_i =r_{\!fi }:= g_i (df   )\in \CO_{X x}$. Then $r_{ i} (x)=0$, 
$\{ dr_{i} (x)\}_{i\le n-1}$ is a base of $\tau^\perp \subset T^*_x X$, and  $dr_{n} (x)=0$. Let $Z=Z_f$ be the (germ of) smooth curve passing through $x$ defined by equations $r_{1} =\ldots = r_{n-1}=0$;  its tangent line $T_x Z  \subset T_x X$ is generated by $\tau$. 
 
One has $\langle C, df   \rangle_{(x,\nu )}=\dim \CO_{X x}/\Sigma \,   \CO_{Xx}r_i = \dim  \CO_{Z x}/  \CO_{Z x} r_n $, so we look for  $f\in\CS$ with $r_{n}|_{Z}$ having zero of order 2 at $x$. If  $f$ we started with does not fit the condition, then we modify it as follows. Pick any  $a\in \CO_{Xx}$  such that $a(x) =0$ and $\tau ( a ) \neq 0$. Then  $f':= f +a^3  $ is what we need. Indeed, $df'=df +a^2 da$, so $f'\in\CS$. Set $r'_i := r_{f' i}$, $Z' :=Z_{f'}$, etc. One has
 $dr'_i (x)=dr_i (x)$ and $r'_n$ equals $ r_n + a^2\tau ( a ) $ modulo the cube of the maximal ideal of $\CO_{X x}$.  Since $Z' $  is tangent to $Z$ at $x$,  $r_{n}|_{Z'}$ has zero of order $>2$ at $x$. So $r'_{n}|_{Z'}$ has zero of order 2; we are done. \hfill$\square$ \enddemo

\subhead{\rm 4.13}\endsubhead {\it Proof of the last assertion of Theorem 1.7.} We are in the setting of 1.7. Let $C$ be any irreducible closed conical subset of $T^* \Bbb P$ of dimension $n$, $D:=D_C$. 

\proclaim{Proposition} The map $\tilde{p}^\vee_C: P(i_\circ C) \to D$
 is   birational if $p\neq 2$.  If $p=2$, then  $\tilde{p}^\vee_C$   can also be of generic degree  2. The latter happens if and only if $C$ is $\kappa$-isotropic\footnote{I.e., $\kappa$ vanishes on $TC_{\text{reg} }\subset T(T^* \Bbb P )$.} and for  $(x,\nu )\in C_{\text{reg}}$ the rank of the   map $ T_{(x,\nu )} C \to T_x \Bbb P $ is odd.\footnote{The parity of the rank does not depend on the choice of $(x,\nu )\in C_{\text{reg}}$:  it is odd if and only if $T_{(x,\nu )}C$ and $T^v (T^* \Bbb P)$ lie in the different components of $\text{\rm Gr}^\kappa (T_{(x,\nu )}(T^* \Bbb P ))$.}
\endproclaim

\remark{Example} If $C=T^*_Y \Bbb P$, then   $ \tilde{p}^\vee_{ C}   $ is {\it not} birational  if and only if $p=2$, $\dim Y$ is odd.
\endremark

\demo{Proof} We know that $\tilde{p}^\vee_C$ is generically radicial. Let $D^o$ be a dense open subset of 
$D_{\text{reg}}$ such that $\tilde{p}^\vee_C$ is flat (hence radicial) over it, and the scheme $P(i_\circ C)^o := P(i_\circ C)_{D^o}$ is smooth. For a closed point $\tilde{x}^\vee \in D^o$ 
  there is a unique    $x\in\Bbb P$ with $(i(x),\tilde{x}^\vee )\in P(i_\circ C)\subset \tilde{Q}$. 
The scheme-theoretic fiber $P(i_\circ C)_{\tilde{x}^\vee}$ is a finite scheme supported at $(i(x),\tilde{x}^\vee )$; its order is the generic degree $\delta$ of $\tilde{p}^\vee_C$.

We use the notation from 4.8. Let  $L$ be a line in $ \tilde{\Bbb P}^\vee$ that intersects $D$ at $\tilde{x}^\vee$ transversally and such that $L\hra \tilde{\Bbb P}^\vee$ is $(i_\circ C)^\vee$-transversal at $\tilde{x}^\vee$. 
Then, by 4.8(iv), one has $x\in U_L$;  choose $t$ and $\tilde{U}$  as in 4.8(ii) so that $\tilde{U}$ does not intersect any other critical fiber $\tilde{Q}_{y}$, $y\in L\cap D$, $y\neq \tilde{x}^\vee$. One has $P(i_\circ C)_{\tilde{x}^\vee}= P(i_\circ C )_{\tilde{U}}\cap \tilde{Q}_L =df_{(L,t)}(U)\cap C$ (the first equality comes since $L$ intersects $D^o$ at $x^\vee$ transversally,   for the second one see 4.8(ii)). Thus $\delta =\langle df_{(L,t)}(U), C\rangle_{(x,\nu )}$ since  $\nu :=df_{(L,t)}(x)\in C_{\text{reg}}$. 

Set $X:=\Bbb P \smallsetminus H$ for $H$ as in 4.8(iii); we view it as a vector space with $x=0$.
Due to 4.8(iii) we are reduced to the next assertion:  For every $\nu \in C_{\text{reg}\, x}\subset T^*_x X$ one can find $f$ as in the proposition in 4.12 such that $f=q_1 /q_2$ where $q_1$, $q_2$ are nonzero non-proportional polynomials of degree $\le d$, $q_2 (x)\neq 0$. To check it look at the proof in 4.12 (we follow the notation there). Since $f$ with $df (x)=\nu$, $d^{(2)}f (T_x X)=\CL$, can be chosen to be a quadratic polynomial $f=\nu +q$ (see the lemma in 4.12), this solves our problem in case (i) of the proposition in 4.12.  In case (ii) of loc.~cit.~we look for $f$ such that $r_{fn}|_{Z_f}$ has zero of order 2 at $x$. We start with  $f=\nu +q$ as above. 
If it does not satisfy the condition, we modify $f$ as follows.
Pick a linear function $a$ such that $\tau (a)\neq 0$. The cubic polynomial $f'=\nu +q+a^3$  is a solution for $d\ge 3$. For $d=2$ we find a solution $f''$ which is the ratio of a quadratic polynomial and a linear one: Namely, 
we can assume that $q|_{Z_f}$ has zero of order 2 at $x$ (otherwise replace $q$ by $q+a^2$); then our $f''$ is $\nu + q (1+a)^{-1}$.   
\hfill$\square$
\enddemo

\end